\documentclass[12pt]{article}
\pdfoutput=1
\usepackage{amsfonts}
\usepackage{amsmath}
\usepackage{amssymb}
\usepackage{graphicx}
\usepackage{epstopdf}
\usepackage{url}
\usepackage[round]{natbib}

\begin{document}

\title{A stable finite-difference scheme for growth and diffusion on a map}

\author{W.~P.~Petersen, S.~Callegari, N.~Tkachenko, \\
J.~D.~Weissmann, Ch.~P.~E.~Zollikofer \\
Anthropological Institute and Museum, University of Z\"{u}rich}

\maketitle

\begin{abstract}
We describe a general Godunov-type splitting for numerical simulations
of the Fisher--Kolmogorov--Petrovski--Piskunov growth and diffusion
equation in two spatial dimensions. In particular, the method is
appropriate for modeling population growth and dispersal on a
terrestrial map. The procedure is semi-implicit, hence quite stable,
and approximately $O(\Delta x^{2}) + O(\Delta t^{2})$ accurate,
excluding boundary condition complications. It also has low memory 
requirements and shows good performance. We illustrate an
application of this solver: global human dispersal in the late
Pleistocene, modeled via growth and diffusion over geographical maps of
paleovegetation and paleoclimate.
\end{abstract}

\begin{flushleft}
\end{flushleft}

\section{Introduction}

There is an increasing interest in modeling population dynamics at
large spatial and temporal scales, for example the modern human
out-of-Africa dispersal \citep{eriksson2012, henn2012, nikitas2005,
  young1995} or Neanderthal dispersal and extinction 
\citep{callegari2013}. 
These models are required to interpret local and global
patterns of genetic, phenetic and cultural variation
\citep{bouckaert2012, eriksson2012, manica2007, pinhasi2005,
  relethford2004}.

Fisher \citep{fisher1937} studied the description -- via a
reaction-diffusion equation -- of an analogous but one-dimensional
problem: the propagation of an advantageous genetic mutation within an
already-present population, situated along a coast line. Kolmogorov,
Petrovskii and Piskunov \citep{kolmogorov1937} were more general; in
particular, their analysis treated the two-dimensional case. Such a
model (called Fisher/KPP in the following) was first applied to the
dispersal and growth of a population by Skellam \citep{skellam1951},
and serves as an important control for designing and
validating other more complex spatiotemporal population 
models \citep{callegari2013}.
Coupling
population dynamics with models of large-scale changes in continental
topography, climate, and ecosystem productivity is essential to
understand the role of environmental constraints on patterns of
genetic, phenetic, and cultural variation among human populations
\citep{callegari2013}.

Here we present a stable and efficient finite-difference solver for
the Fisher/ KPP equation on 2-D domains of arbitrary shape
(e.g.~geographical maps), and show how it can be extended to include
environmental fluctuations. In a brief outline of our paper, we will: 
review the derivation of the Fisher/KPP equation (Section~\ref{sec:fisherkppeq}); 
develop finite-difference schemes in 1 and
2 dimensions for constant environmental carrying capacity, $\mathcal{K}$
(Section~\ref{sec:methodsandsplittings}); extend the scheme to allow for space- and
time-dependent $\mathcal{K}(\mathbf{x},t)$ (Section~\ref{sec:fisherkpp})
and irregular domains such as geographical maps; 
show an application of this technique to the out-of-Africa
dispersal of \emph{Homo sapiens} by using net primal productivity
(NPP) as a proxy for $\mathcal{K}(\mathbf{x},t)$ (Section~\ref{sec:dispersal}).

\section{The Fisher/KPP equation}
\label{sec:fisherkppeq}
An intuitive way to get the Fisher/KPP equation
\citep{fisher1937,kolmogorov1937,skellam1951} is as follows. A current
$\mathbf{j}$ of particles (e.g., individuals) moving across an
interface located at $\mathbf{x}$ is proportional to the gradient of
the population density $p$ (``Fickian diffusion'')
\[ \mathbf{j} = -c\,\nabla p. \]
%Without any source of new population, 
The rate of change of $p$ is then given by the mass balance equation
\citep{reichl2004},
%\[ \frac{\partial p}{\partial t} + \nabla \cdot \mathbf{j} = \rho, \]
which for Fickian diffusion reads
%Thus, if there is no source term, we have
%\[ \frac{\partial p}{\partial t} - c\,\nabla^2 p = 0, $$ which is the
%heat equation when $c = \mathcal D/2$ and $\mathcal D$ is the
%diffusion coefficient. For lack of a better model, we assume $c$ is a
%constant (Young \& Bettinger).  We now expand the notion of Fickian
%diffusion to include a source:
\[ 
\frac{\partial p}{\partial t} - c\,\nabla^2 p = \rho.
\] 
If $\rho=0$, this is the heat equation when $c = \mathcal D/2$ 
and $\mathcal D$ is the
diffusion coefficient. For lack of a better model, we assume $c$ is a
constant (Young \& Bettinger \citep{young1995}). 
The source term $\rho$ is usually
modeled by a logistic growth function, $\rho = \lambda p (1-p/\mathcal
K)$, yielding the Fisher/KPP equation
\begin{equation}
\frac{\partial p}{\partial t} = \lambda
\left(1-\frac{p}{\mathcal{K}}\right) p + c\,\nabla^2 p,
\label{eq:fisherKPP0}
\end{equation}
where $\mathcal{K}$ is called \emph{carrying capacity} and $\lambda$
is the growth rate. In (\ref{eq:fisherKPP0}) we assumed $\mathcal K$ is
constant, but it suffices that there is an \emph{upper limit}
$\bar{\mathcal K}$, in which case $0 \leq p / \bar{\mathcal K} \leq
1$. This scaled version will be used in (\ref{eq:fisherKPP}) below.
%Space coordinates $\mathbf{x}$ can be scaled so that $c = \mathcal D/2
%=1/2$; these scalings are given in Table~\ref{tab:scalings}.  Both
Fisher and KPP were particularly interested in the traveling wave case,
$p(t,\mathbf{x}) = f(\mathbf{x}_{0} + \mathbf{v} t)$. Notice what
happens here if $f$ exists:
\[
\mathbf{v} \cdot \nabla f = \lambda
\left(1-\frac{f}{\mathcal{K}}\right) f + c\,\nabla^2 f,
\]
which in 1-D
%, by rescaling $x \leftarrow x \cdot \sqrt{2 c}$ to 
%set $c=1/2$, 
becomes a second order ordinary differential equation 
\[
v \frac{{\rm d} f}{{\rm d} x} = \lambda
\left(1-\frac{f}{\mathcal{K}}\right) f + c \frac{{\rm d}^2
  f}{{\rm d} x^2}.
\]
If $v$ were known, this ODE could be solved using
\texttt{pvp4c} from \texttt{MatLab}, for example. 
%By rescaling 
%$p(t,\mathbf{x}) = \mathcal{K} u(t,\mathbf{x})$, 
%equation (\ref{eq:fisherKPP0}) becomes
%\[
%\frac{\partial u}{\partial t} = \lambda (1-u)u + \frac{1}{2}\,\nabla^2
%u
%\]
%Finally, rescaling $t \leftarrow \lambda t$ and $\mathbf{x} \leftarrow
%\sqrt\lambda \mathbf{x}$, the parameter $\lambda$ can be set to
%unity. Hence, 
By the rescalings show in Table~\ref{tab:scalings}, for constant
$\mathcal{K}$ the Fisher/KPP equation (\ref{eq:fisherKPP0}) will be
used in the form
\begin{equation}
\frac{\partial u}{\partial t} = (1-u)u + \frac{1}{2}\,\nabla^2 u,
\label{eq:fisherKPP}
\end{equation}
where the only sensible solutions have $0 \leq u \leq 1$.  The initial
distribution $u(0,\mathbf{x}) = u_{0}(\mathbf{x})$ must be defined for
all $\mathbf{x}$. 

\begin{table}[t]
\caption{Variables in left column are scaled versions of those in
  right column.\label{tab:scalings}}
\begin{center} \footnotesize
\renewcommand{\arraystretch}{1.3}
\begin{tabular}{|c|c|} \hline
In (\ref{eq:fisherKPP}) & In (\ref{eq:fisherKPP0}) \\ \hline $
\mathbf{x}$ & $\sqrt{\frac{\lambda}{2 c}} \, \mathbf{x}$ \\ $t$ &
$\lambda t$ \\ $u$ & $p/\mathcal{K}$ \\ \hline
\end{tabular}
\end{center}
\end{table}

In Murray (\citep{murray2002}, eq.~(11.17)) our growth coefficient
$\lambda$ is called $r$ and $c$ is denoted by $D$, whereas in Young
and Bettinger \citep{young1995} the growth coefficient is $R$ and the
diffusion coefficient is $K$. These inputs to our code are given in
units of yr$^{-1}$ and km$^2$/yr respectively.

\section{Numerical methods and splitting}
\label{sec:methodsandsplittings}
In one dimension, (\ref{eq:fisherKPP}) can be solved using the
\texttt{MatLab} function \texttt{pdepe}. 
In fact, if the system is two-dimensional
but rotationally symmetric, \texttt{pdepe} can again be used with the
radial part of the Laplace operator in cylindrical coordinates,
\[
\nabla^2 = \frac{1}{r} \frac{\partial}{\partial r} r \frac{\partial}{\partial r}
+ \mbox{non-contributing terms},
\]
requiring only that one sets a \texttt{pdepe} parameter \texttt{m}=1.
Although the \texttt{MatLab} function \texttt{pdepe} is robust, it cannot be
generalized to arbitrary 2-D domains. However, it is a
valuable control for testing more general solvers, and a more general
solver is what we wish to explore here. 

\subsection{The finite-difference scheme}
Since the map on which we will be working is a pixelized plane, an obvious
method uses finite differences. 
First, however, let us examine the 1-D case for (\ref{eq:fisherKPP}).  
In this situation, the
second order derivative becomes a differencing operator in matrix form
acting on the vector $\{ u_{j}, j=1,n \}$, where $u_{j} =
u(x_{0}+(j-1) \Delta x)$,
\[
\frac{d^{2} u}{dx^{2}} \rightarrow \frac{1}{(\Delta x)^{2}} A u,
\]
where the matrix $A$ is
\begin{equation}
A = \left( \begin{array}{cccccc}
-2 & 1 & 0 & 0 & \ldots & 0 \\
1 & -2 & 1 & 0 & \ldots & 0 \\
0 & 1 & -2 & 1 & \ldots &  \\
0 & 0 & 1 & \ldots &  &  \\ 
0 & 0 & \ldots & 0 & 1 & -2 \end{array} \right). \label{eq:Amatrix}
\end{equation}
If $h$ is the time step, the Courant–-Friedrichs–-Lewy (CFL) parameter
\citep{strang1986} is
\[
k = \frac{h}{2 (\Delta x)^{2}}.
\] An explicit integrator for (\ref{eq:fisherKPP}) would require $k <
1/4$ \citep{leveque2007, strang1986}.  In our case, because the
boundary conditions are so irregular on a map, we are less interested
in a method of higher order than $2nd$ because stability is more important
\citep{godunov1987}.

Using this notation, the lowest order approximation is Euler's method
which estimates the next step $u(t+h)$ by
\begin{equation}
u_{E} = u(t) + k \, A u(t) + h \, (1-u(t)) u(t),
\label{eq:Eulerestimate}
\end{equation}
which should be considered a vector equation in $u(t) = \{ u_{j}(t),
j=1,n \}$. The logistic terms, which are diagonal, should be taken to
mean $((1-u)u)_{j} = (1-u_{j})u_{j}$ for $j=1,2, \ldots, n$.  Euler's
method is both low-accuracy and usually unstable if it is used alone
over many steps.  But, it is $O(h)$ accurate and thus useful as an
explicit estimate in $O(h)$ terms.  An application of the
trapezoidal rule yields
\begin{eqnarray}
u(t+h) & = & u(t) + \frac{k}{2} \left( A u(t+h) + A u(t)
\right) \label{eq:trapezoidalrule} \\ & & \qquad + \frac{h}{2}
\left((1-u(t+h))u(t+h) + (1-u(t))u(t) \right) \nonumber
\end{eqnarray}
and is an $O(h^{2}) + O((\Delta x)^{2})$ accurate procedure but
solving the quadratic vector equation (\ref{eq:trapezoidalrule}) for
$u(t+h)$ is awkward. To the same $O(h^{2})$, we propose a
semi-implicit procedure which uses the Euler estimate
(\ref{eq:Eulerestimate}) to replace one of the terms in
(\ref{eq:trapezoidalrule}):
\begin{eqnarray}
u(t+h) & = & u(t) + \frac{k}{2} \left( A u(t+h) + A u(t) \right) \label{eq:semiimplicit}\\
& & \qquad + \frac{h}{2} \left( (1-u_{E})u(t+h) + (1-u(t))u(t) \right).\nonumber
\end{eqnarray}
Equation (\ref{eq:semiimplicit}) can now be solved as a linear system,
\begin{equation}
\left( 1 - \frac{k}{2} A - \frac{h}{2} (1-u_{E}) \right) u(t+h) = 
u(t) + \frac{k}{2} A u(t) + \frac{h}{2} (1-u(t))u(t), 
\label{eq:semidifference}
\end{equation}
because the matrix, $1 - \frac{k}{2} A - \frac{h}{2} (1-u_{E})$,
on the left hand side is explicit, as is the right hand side.
That is, this matrix and the right hand side contain only old data,
namely only information from the previous step, $u(t)$. Euler estimate
$u_{E}$ is an explicit one step computation using $u(t)$.
Significant advantages are: the matrix on the left hand side
is tridiagonal with constants on the sub/super-diagonals, 
and the diagonal terms are $O(1)$ strong. 
The procedure (\ref{eq:semidifference}) is only linearly stable
but we will show empirically that it gives good results
when compared to \texttt{pdepe} when this \texttt{MatLab} function is
appropriate, that is, in both the one-dimensional and rotationally
symmetric 2-D case. Not only is the method (\ref{eq:semidifference})
step-wise stable but also stable for initial data which may not be smooth.

\begin{figure}
\begin{center}
\begin{minipage}{0.48\textwidth}
  \includegraphics[width=\textwidth]{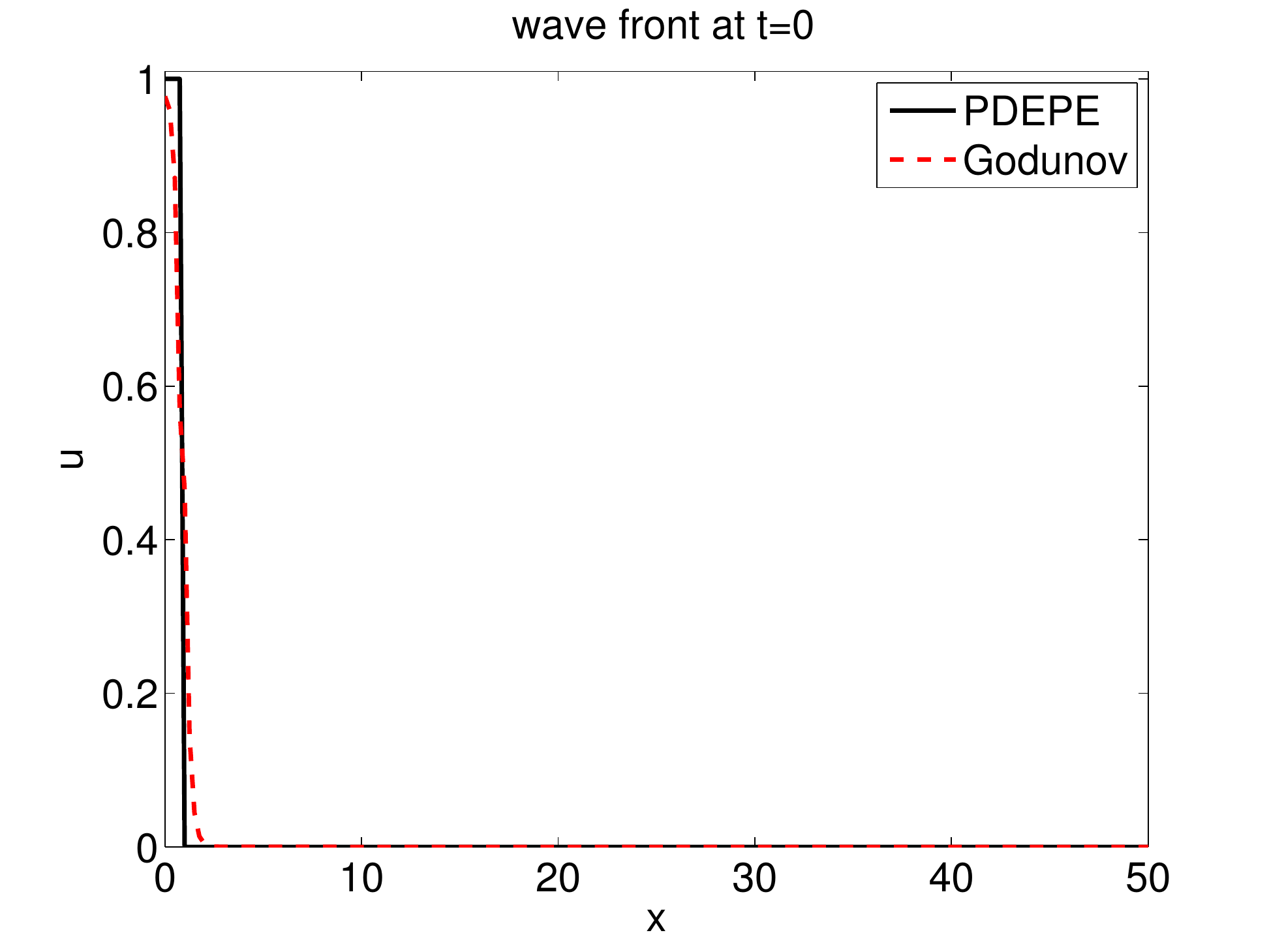}
\end{minipage}
\hfill
\begin{minipage}{0.48\textwidth}
  \includegraphics[width=\textwidth]{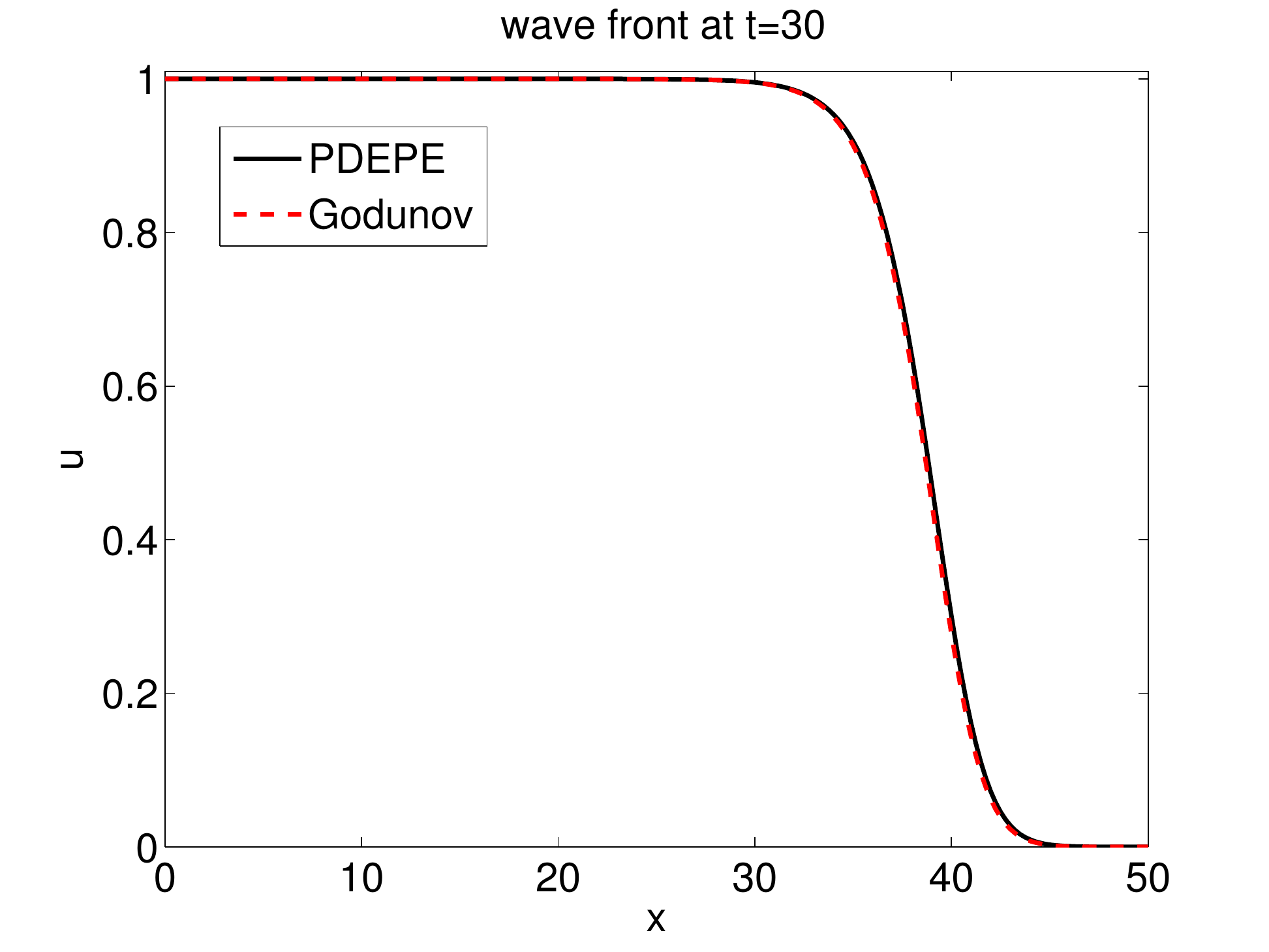}
\end{minipage}
  \caption{Left: Godunov vs. \texttt{pdepe} at the end 
   of one time step $t=h$. Right: Godunov vs. \texttt{pedpe} 
   at time $t=30$. The time step $h=1/5$ and $\Delta x = 1/5$.}
\label{fig:1dplots}
\end{center}
\end{figure}

Figure~\ref{fig:1dplots} shows the results for $h=1/5$, $k=2.5$
compared to \texttt{pdepe}. Notice that at $t=20$ the agreement
is remarkable; and that at $t=h$, where the wave front profile is very steep,
our Godunov splitting described in 
Section~\ref{sec:Godunov2Dsplitting}, specifically
eq.(\ref{eq:yoshida}), is very stable. The CFL number, $k=2.5$,
used to get Figure~\ref{fig:1dplots} is much larger
than would be possible with an explicit method
\citep{leveque2007}.

\subsection{2-D case: Godunov--Strang--Yoshida splittings}
\label{sec:Godunov2Dsplitting}
It turns out that a generalization to the 2-D problem is a
straightforward variant of Strang--Yoshida splittings
\citep{strang1968, yoshida1990}, which are themselves variants of
Godunov's method \citep{godunov1987}.  The following is a fully
implicit variant of our two-dimensional scheme, with two intermediate
arrays, $u^{\star}$ and $u^{\star \star}$,
\begin{subequations} \label{eq:yoshida}
\begin{align}
u^{\star} & = u(t) + \frac{k}{4} \left( A_{x} u^{\star} + A_{x} u(t) \right)
\label{eq:yoshidaA} \\ 
u^{\star \star} & = u^{\star} + \frac{k}{2} \left( A_{y} u^{\star \star} + A_{y} u^{\star} \right)
\label{eq:yoshidaB} \\
& \qquad \; + \frac{h}{2} \left( (1-u^{\star \star}) u^{\star \star} + (1-u^{\star}) u^{\star} \right)
\nonumber \\ 
u(t+h) & = u^{\star \star} + \frac{k}{4} \left( A_{x} u(t+h) + A_{x} u^{\star \star} \right).
\label{eq:yoshidaC} 
\end{align} 
\end{subequations}

In (\ref{eq:yoshida}), the operators $A_{x}$ and $A_{y}$ are the same
as (\ref{eq:Amatrix}) for directions $x$ and $y$, respectively. For
simulations on a lattice, $u_{ij}(t) = u(t,x_{0}+(i-1)\Delta
x,y_{0}+(j-1)\Delta y)$, where $1 \leq i \leq N_x,\, 1 \leq j \leq
N_y$ and $\Delta x = \Delta y$, the following gives the action of the
$A_{x}, A_{y}$ operators:
\begin{eqnarray*}
A_{x} u_{i,j} & = & u_{i-1,j} - 2 u_{i,j} + u_{i+1,j}, \\
A_{y} u_{i,j} & = & u_{i,j-1} - 2 u_{i,j} + u_{i,j+1}.
\end{eqnarray*}
The compression scheme and code outline given
in Appendix~\ref{sec:code} show that only a maximum
of one row or column (i.e., max($N_x,N_y)$) of storage is needed for
$u^{\star}$ and $u^{\star \star}$.

Again because the fully implicit quadratic vector equation in
(\ref{eq:yoshidaB}) is awkward to solve, we use an Euler estimate in
one of the terms. Here is one integration time step of
(\ref{eq:yoshida}) in discrete semi-implicit form:
\begin{subequations} \label{eq:discreteyoshida}
\begin{align} 
\left( 1 - \frac{k}{4} A_{x} \right) u^{\star} & = \left( 1 +
\frac{k}{4} A_{x} \right) u(t) \label{eq:discreteyoshidaA} \\ u_{E} &=
u^{\star} + k A_{y} u^{\star} + h (1-u^{\star})
u^{\star} \label{eq:discreteyoshidaB} \\ \left( 1 - \frac{k}{2} A_{y}
- \frac{h}{2} (1-u_{E}) \right) u^{\star \star} & = \left( 1 +
\frac{k}{2} A_{y} + \frac{h}{2} (1-u^{\star}) \right) u^{\star}
\label{eq:discreteyoshidaC} \\
\left( 1 - \frac{k}{4} A_{x} \right) u(t+h) & = \left( 1 + \frac{k}{4}
A_{x} \right) u^{\star \star}.
\label{eq:discreteyoshidaD}
\end{align}
\end{subequations}
Equations (\ref{eq:discreteyoshidaA}), (\ref{eq:discreteyoshidaC}),
and (\ref{eq:discreteyoshidaD}) are solved in sequence as multiple
independent tridiagonal systems for $u^{\star}, u^{\star \star}$ and
the final step $u(t+h)$.

\subsection{Symmetries in 2-D case}
\begin{figure}
\begin{center}
\begin{minipage}{0.48\textwidth}
  \includegraphics[width=\textwidth]{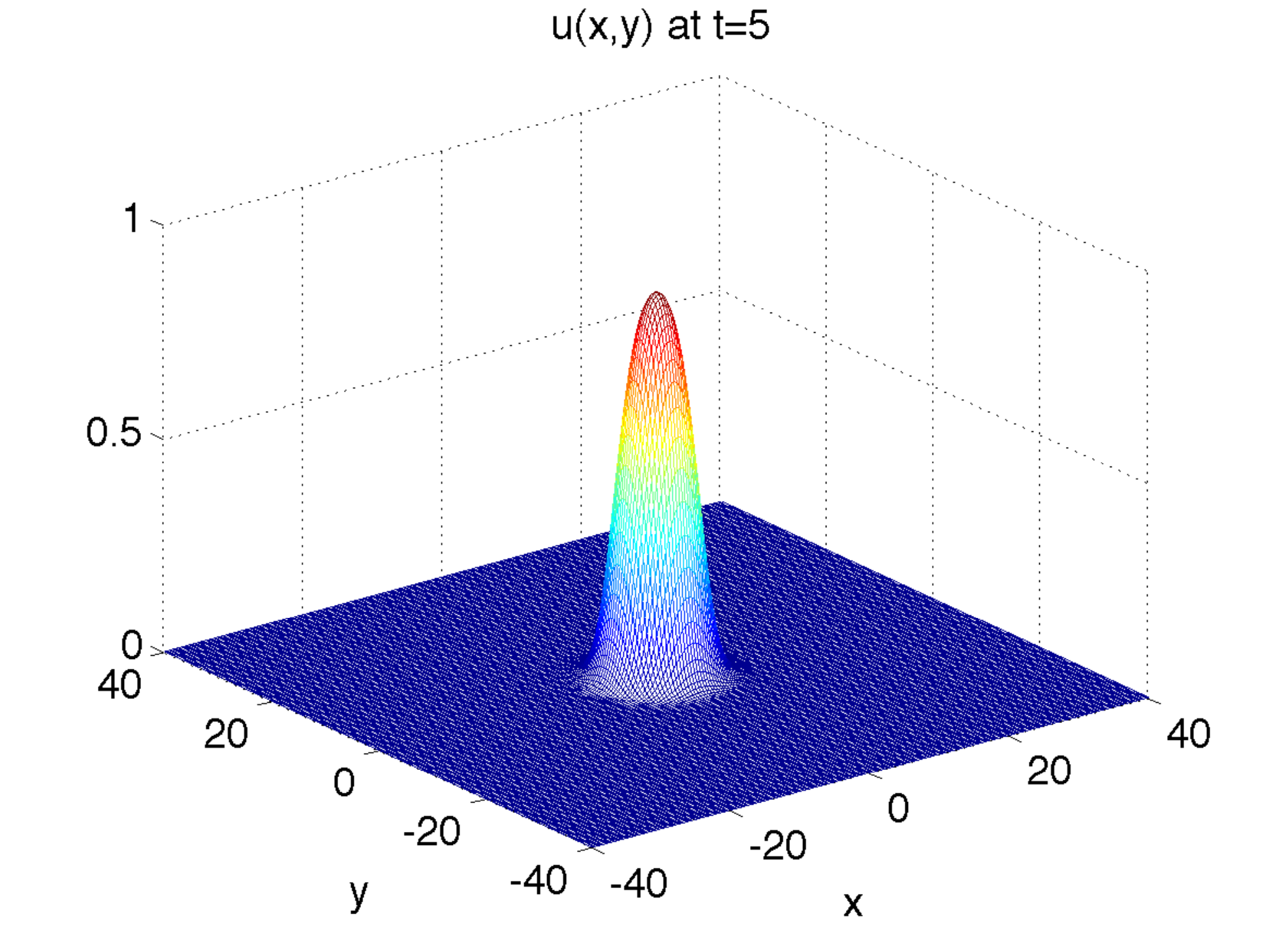}
\end{minipage}
\hfill
\begin{minipage}{0.48\textwidth}
  \includegraphics[width=\textwidth]{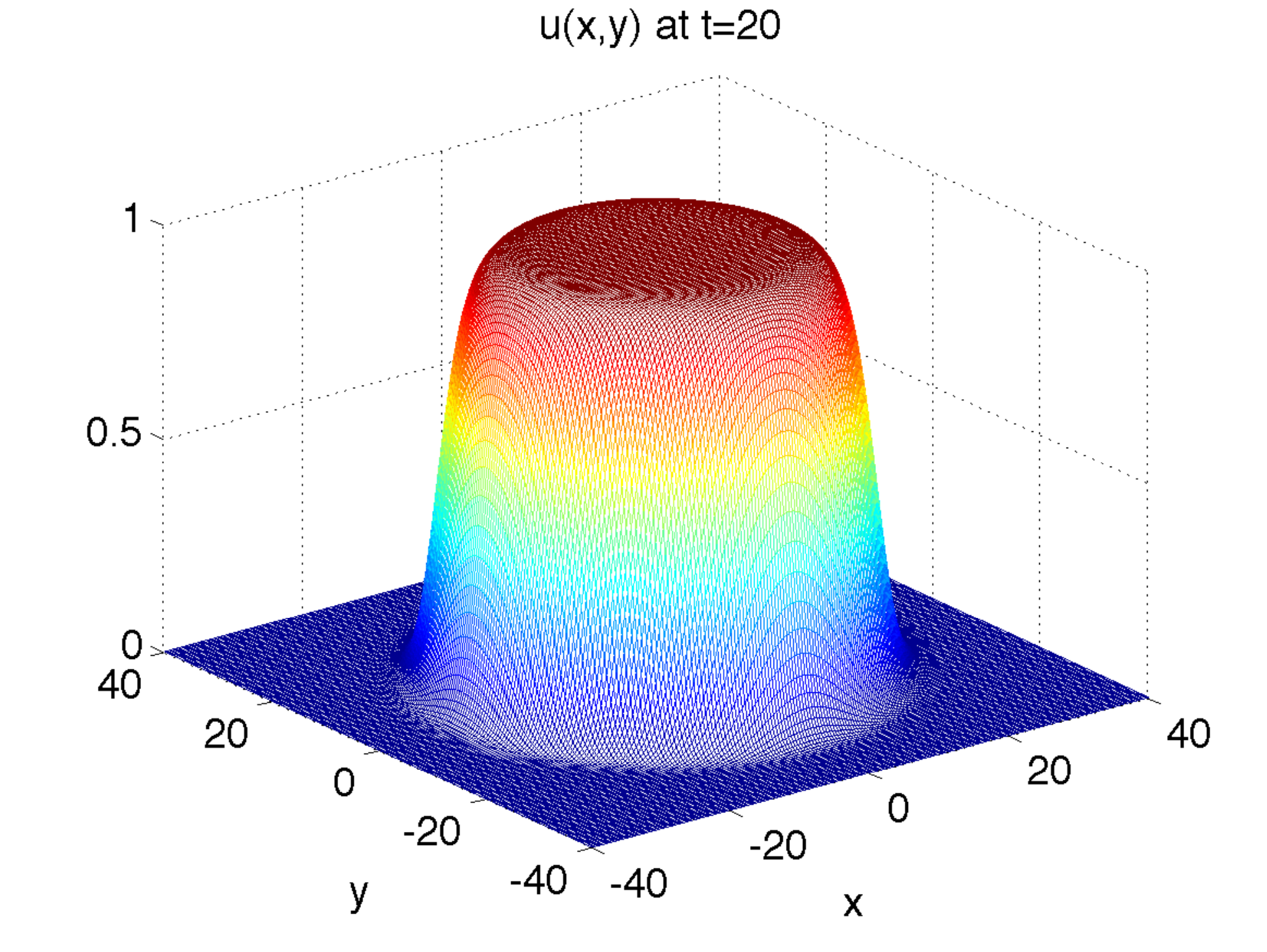}
\end{minipage}
\caption{Left: 2-D solution at $t=5$ obtained with Godunov scheme
  (\ref{eq:discreteyoshida}). Right: same, at $t=20$. The time step
  was $h=1/10$, $\Delta x = 2/5$ and $N_X=N_Y=201$.}
\label{fig:2dpcolorplots}
\end{center}
\end{figure}

\begin{figure}
\begin{center}
\begin{minipage}{0.48\textwidth}
  \includegraphics[width=\textwidth]{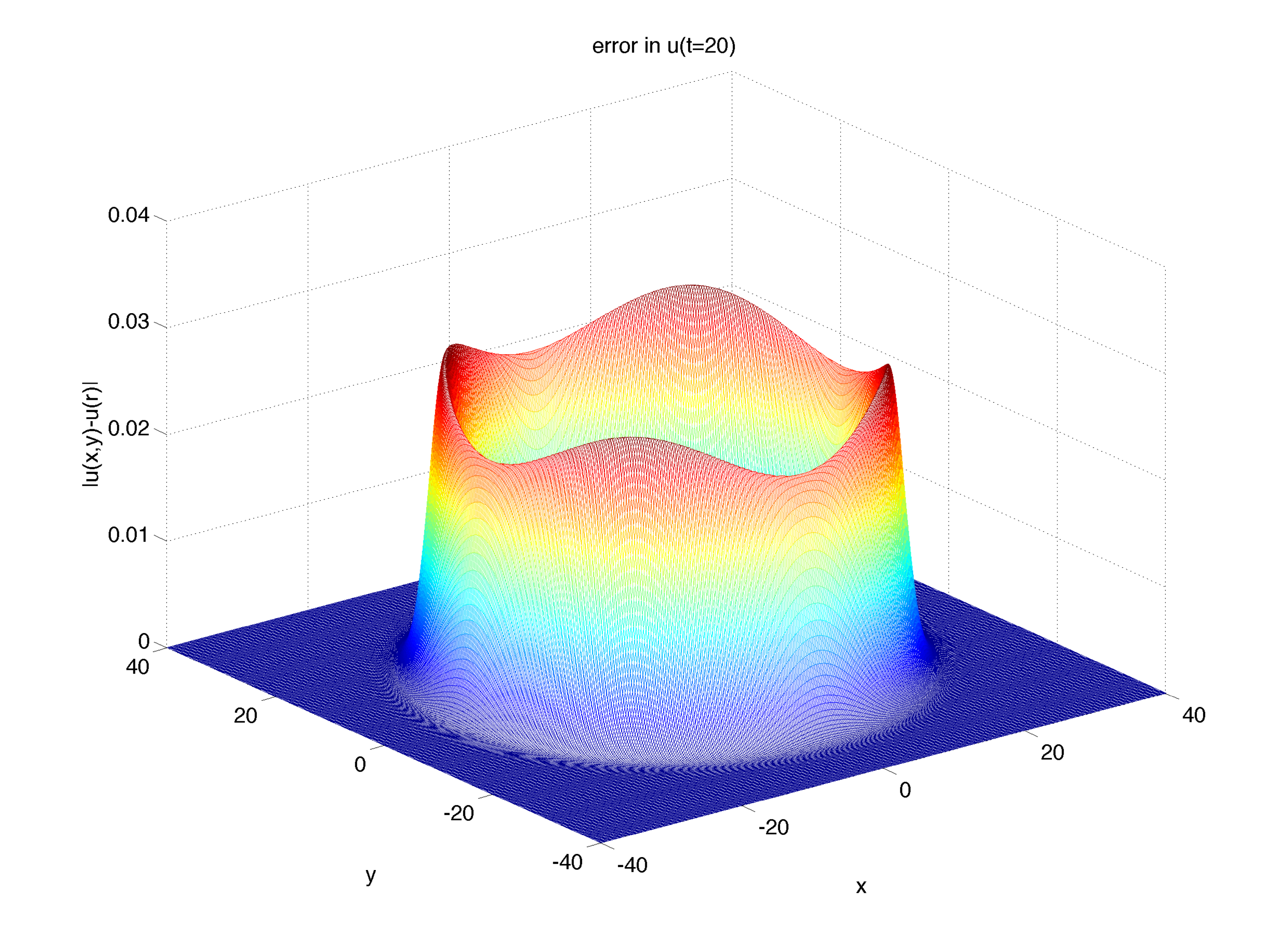}
\end{minipage}
\hfill
\begin{minipage}{0.48\textwidth}
  \includegraphics[width=\textwidth]{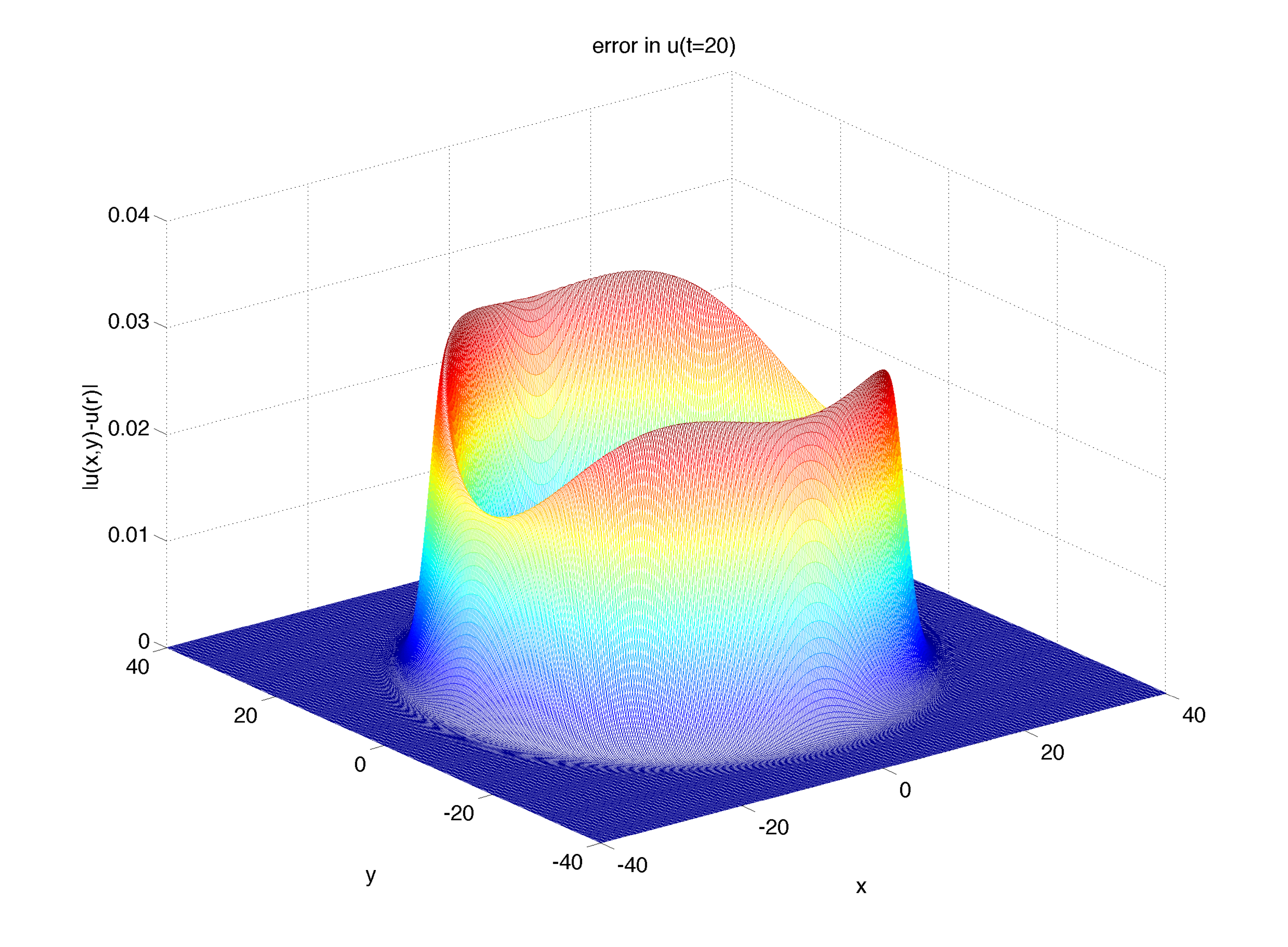}
\end{minipage}
\caption{Left: $|u(x,y)-u(r)|$ error profile of the Godunov scheme
  compared to the \texttt{pdepe} solution at $t=20$ symmetrizing $x
  \leftrightarrow y$ on alternate time steps: $\epsilon_{RMS} = 7.3
  \cdot 10^{-3}$, $\epsilon_{MAX} = 2.9 \cdot 10^{-2}$.  Right: the
  same, but with fixed $x$ and $y$: $\epsilon_{RMS} = 7.4 \cdot
  10^{-3}$, $\epsilon_{MAX} = 3.0 \cdot 10^{-2}$.  Step size is
  $h=1/10$, $N_x=N_y=401$.}
\label{fig:symmetryornot}
\end{center}
\end{figure}

Our Godunov scheme (\ref{eq:discreteyoshida}) is not rotationally
symmetric, and thus one way to estimate the error is to assess a
solution using (\ref{eq:discreteyoshida}) for a symmetric
problem. Again, we can use \texttt{pdepe} but now with the
cylindrically symmetric parameter choice \texttt{m}=1 (see Section
12.5 in \citep{higham2005}). Figure~\ref{fig:2dpcolorplots} shows
that any asymmetries are not apparent without more careful
examination. Even the wave front portrait of the 2-D case in the
left-hand panel of Figure~\ref{fig:2Dwavefrontanderrors} and the error
estimate in the right-hand panel of the same Figure are not
sufficiently quantitative. In particular, there should be no
distinction between $x$ and $y$ directions in (\ref{eq:yoshida}),
while a 2-D plot of the error distribution shows a small asymmetry
(compare the right-hand plot in Figure~\ref{fig:symmetryornot} to the left). 

For this reason, we implemented an alternating direction method, {\em
  \`a la} Crank-Nicholson \citep{leveque2007,ritchmyer1967}, which
makes the error distribution more symmetric. The left panel of
Figure~\ref{fig:symmetryornot} shows that, while symmetrization only
slightly improves the r.m.s.~and maximum errors, they have now the
desirable property of being more rotationally equi-distributed:
respectively, cyclic groups $C_{4}$ vs. $C_{2}$. Relatively larger
deviations from the \texttt{pdepe} solution now correspond to
directions diagonal to the spatial lattice, as expected, and do not
reflect the arbitrary choice of $x$ and $y$ in the integration.

\section{Fisher/KPP on maps}
\label{sec:fisherkpp}
The next natural step when applying a reaction-diffusion equation to
the modeling of population dispersal is to include geographical and
environmental effects. In this Section, we discuss how to implement
our solver on domains with space- and time-dependent $\mathcal K$, and
then how to treat irregular boundaries that arise when solving
Fisher/KPP on geographical maps. Using the same Godunov-type splitting
described above, it is more straightforward to do the simulations on a
map than might be expected.

\subsection{Maps with space-dependent capacity}
\label{sec:fisherkppvarK}
In our scaling of (\ref{eq:fisherKPP}), the maximum population density
at $\mathbf{x}$ is unity. Thus, in the following we will use a scaled
carrying capacity $0 < \mathcal{K}(\mathbf{x}) \leq 1$. Our algorithm
(\ref{eq:yoshida}) can be modified in a straightforward way for the
case that $\mathcal{K}(\mathbf{x},t)$ also depends explicitly on time:
see Section~\ref{sec:fisherkppvarKxt}. First, however, let us deal
with the time independent $\mathcal{K}(\mathbf{x})$ case,
\begin{equation}
\frac{\partial u}{\partial t} = (1-\frac{u}{\mathcal{K}}) u +
\frac{1}{2}\,\nabla^2 u,
\label{eq:fisherwithK}
\end{equation}
(we dropped the $\mathbf x$ dependence of $\mathcal K$ for convenience
of notation), for which the Godunov splitting (\ref{eq:yoshida}) now
becomes
\begin{subequations} \label{eq:godunovwithK}
\begin{align}
\left( 1 - \frac{k}{4} A_{x} \right) u^{\star} & = \left( 1 + \frac{k}{4} A_{x} \right) u(t) \label{eq:godunovwithKA} \\
u_{E} &= u^{\star} + k A_{y} u^{\star} + h (1-\frac{u^{\star}}{\mathcal{K}}) 
u^{\star} \label{eq:godunovwithKB} \\
\left( 1 - \frac{k}{2} A_{y} - \frac{h}{2} (1-\frac{u_{E}}{\mathcal{K}}) \right) u^{\star \star} & =
\left( 1 + \frac{k}{2} A_{y} + \frac{h}{2} (1-\frac{u^{\star}}{\mathcal{K}}) \right) u^{\star}
\label{eq:godunovwithKC} \\
\left( 1 - \frac{k}{4} A_{x} \right) u(t+h) & = \left( 1 + \frac{k}{4}
A_{x} \right) u^{\star \star}.
\label{eq:godunovwithKD}
\end{align}
\end{subequations}
Again, as in (\ref{eq:yoshida}), multiple tridiagonal system solves
must be carried out: $N_y$ $x$-direction solutions
(\ref{eq:godunovwithKA}), $N_x$ $y$-direction solutions
(\ref{eq:godunovwithKC}), and finally, another $N_y$
$x$-direction solutions (\ref{eq:godunovwithKD}).  Hopefully no
confusion will result from the notation: $k$ is the CFL parameter,
while $\mathcal{K}$ is the (space-dependent) carrying capacity.

\subsection{A desert test of space-dependent capacity}
\label{sec:gaptest}
Now we are in uncharted territory. To assess if the solver
(\ref{eq:godunovwithK}) works, let us examine a problem where we can
compute a solution by independent means. The test setup follows below.
Its motivations will be explained further in Section~\ref{sec:noisyNPPmaps}.

\noindent
The desert test:
\begin{enumerate}
\item for $-y_{0} \leq y \leq y_{0}$,
  $\mathcal{K}=\mathcal{K}(x)$ is independent of $y$;
\item for $x_{0} \leq x \leq x_{L}$, let $\mathcal{K}=1$, while for
  $x_{L} < x < x_{H}$, let $\mathcal{K}=f_{r}$, where the fraction
  $0 < f_{r} \leq 1$ defines a {\em desert} (inhospitable region) in
  the domain.  Finally, for $x_{H} \leq x \leq x_{1}$, again set
  $\mathcal{K}=1$;
\item initialize $u(x,y,t=0)$ to a strip midway between $x_{0}$ and $x_{L}$,
then run the simulation to study the traveling wave behavior across the $[ x_{L}, x_{H} ]$
desert. 
\end{enumerate}

\begin{figure}
\begin{center}
\begin{minipage}{0.48\textwidth}
  \includegraphics[width=\textwidth]{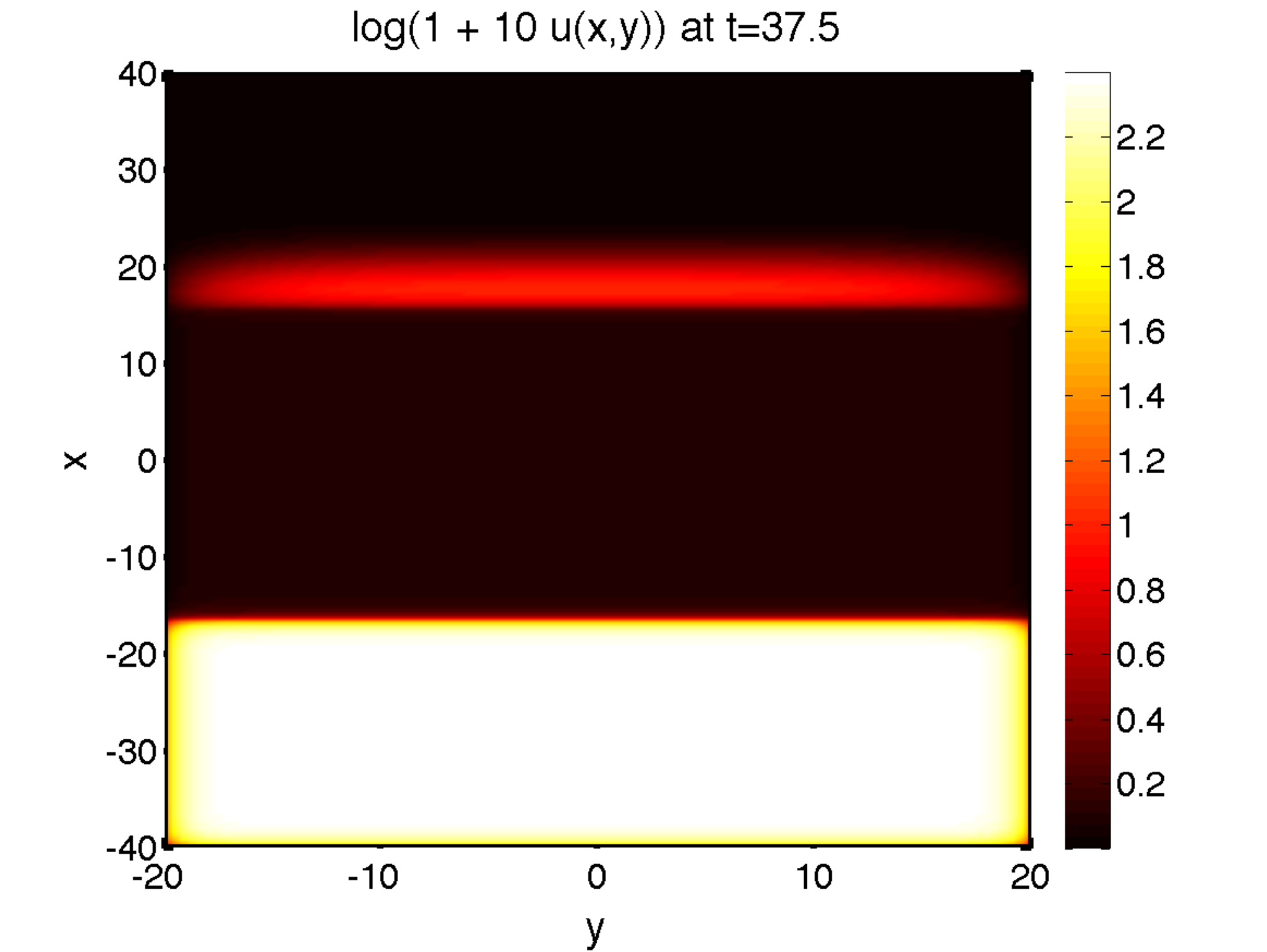}
\end{minipage}
\hfill
\begin{minipage}{0.48\textwidth}
  \includegraphics[width=\textwidth]{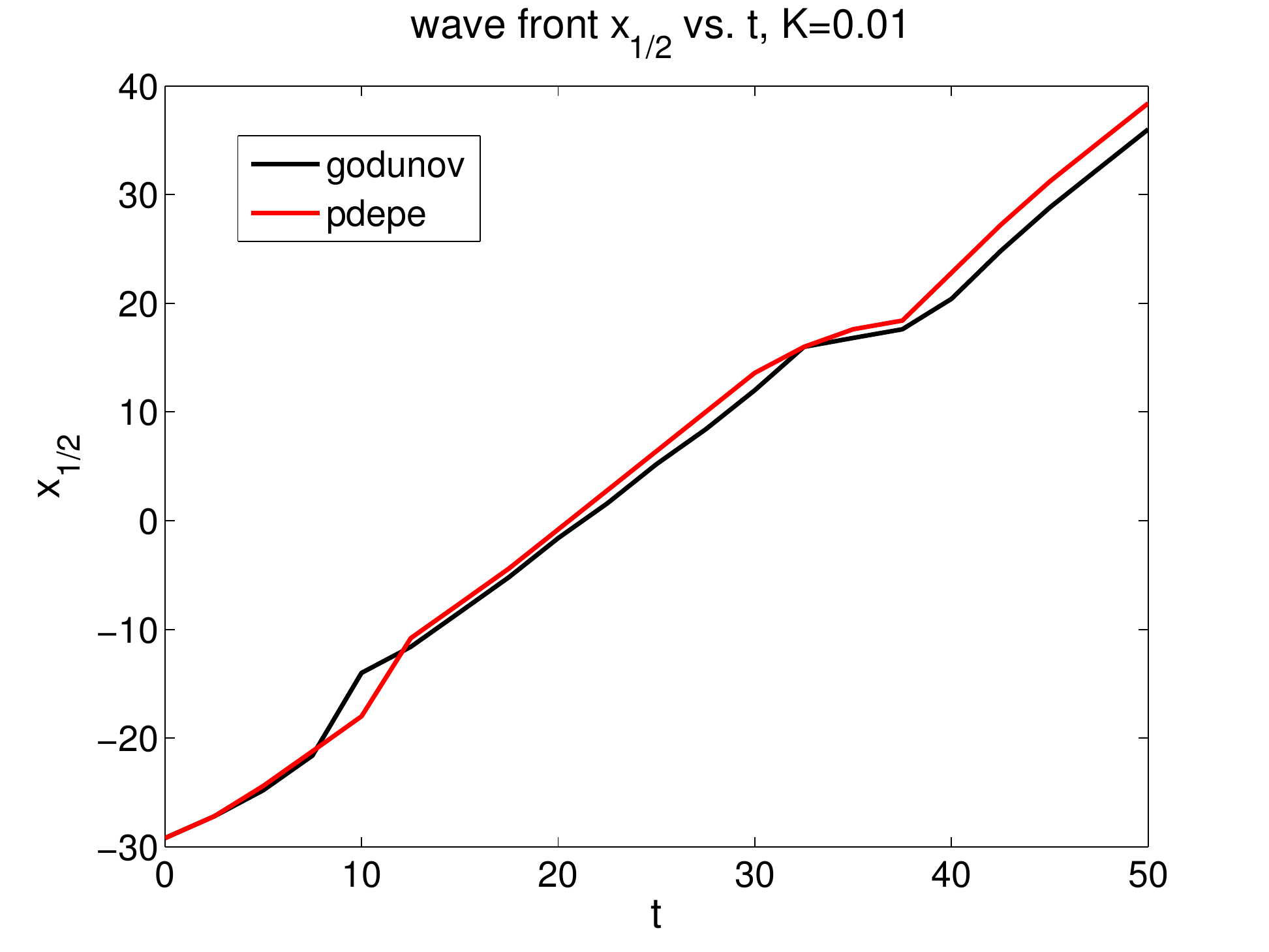}
\end{minipage}
\caption{Left: desert test density at $t=37.5$: $N_x=201$, $N_y=101$,
  and $h=1/8$. For $x \leq -16$ and $x \geq 16$, $\mathcal{K}=1$,
  while $\mathcal{K}=f_{r} = 0.01$ in the $|x| < 16$ desert. The
  image shows $u(x,y)$ shaded from 0 (black) to 1 (white). The wave
  can be seen crossing into the upper $\mathcal{K}=1$ region, while
  the population in the lower $\mathcal{K}=1$ region has already
  saturated. Right: half-height $x_{1/2}(t)$ of the wave front
  vs. time, compared with \texttt{pdepe}.  Velocity $V$ has a
  transient increase, then decrease, as the wave enters/exits the
  desert. Initial condition: a strip with $\langle (x-x_{I})^{2}
  \rangle ^{1/2} =3$ initial population starting at $x_{I}=-30$.}
\label{fig:gaptestresults}
\end{center}
\end{figure}

In other words, the variable carrying capacity domain has
$\mathcal{K}=1$ for $x \leq x_{L}$ and $x  \geq x_{H}$, but
$\mathcal{K}=f_{r}$ in an $x$-direction desert. If the initial data
$u(\mathbf{x},0)$ are widely distributed enough (nearly full
$y$-width), near the middle of the domain, i.e. $y = 0$, the problem
looks basically one-dimensional. Thus we can again use \texttt{pdepe}
from \texttt{MatLab} to compute the behavior of the one-directional
wave as it passes through the desert, and compare this solution to the
behavior of our Godunov method near this same center line. According
to KPP \citep{kolmogorov1937}, the velocity in the desert is the same as
in the ${\mathcal K}=1$ region, and far enough from transients it
should be approximately
\[
V = \sqrt{2 a},
\] where parameter $a=\frac{{\rm d} \rho}{{\rm d} u}(u=0)$. In our
case $\rho(u) = u (1 - u/\mathcal{K})$, so $\frac{{\rm d}\rho}{\rm d
  u}(0) = 1$. Thus, when not entering or leaving the desert, the
velocity should be $\sqrt{2}$ \citep{kolmogorov1937}.
Figure~\ref{fig:gaptestresults} shows that the velocity agrees with
the KPP prediction, and is nearly constant except for short transients
entering and exiting the desert, as expected. Our metric for measuring
this velocity is to find, on the leading edge of the wave front, the
position $x_{1/2}$ where $u(x_{1/2})=\mathcal K/2$.

If the jump in $\mathcal K$ is deep enough ($f_r\ll1$), however, the
integrator will fail without some regularization. The parameters for
the results shown in Figure~\ref{fig:gaptestresults} are
$\mathcal{K}=1$ above and below the desert, but $\mathcal{K}=0.01$ in
the desert. For this case, $\frac{u}{\mathcal{K}} \gg \frac{2}{h}$
gives an instability with unpleasant sign changes, and a
regularization scheme has to be used. We show this in the next
subsection.

\subsection{Regularization against holes and noise}
\label{sec:noisyNPPmaps}
\begin{figure}
\begin{center}
  \includegraphics[width=0.7\textwidth]{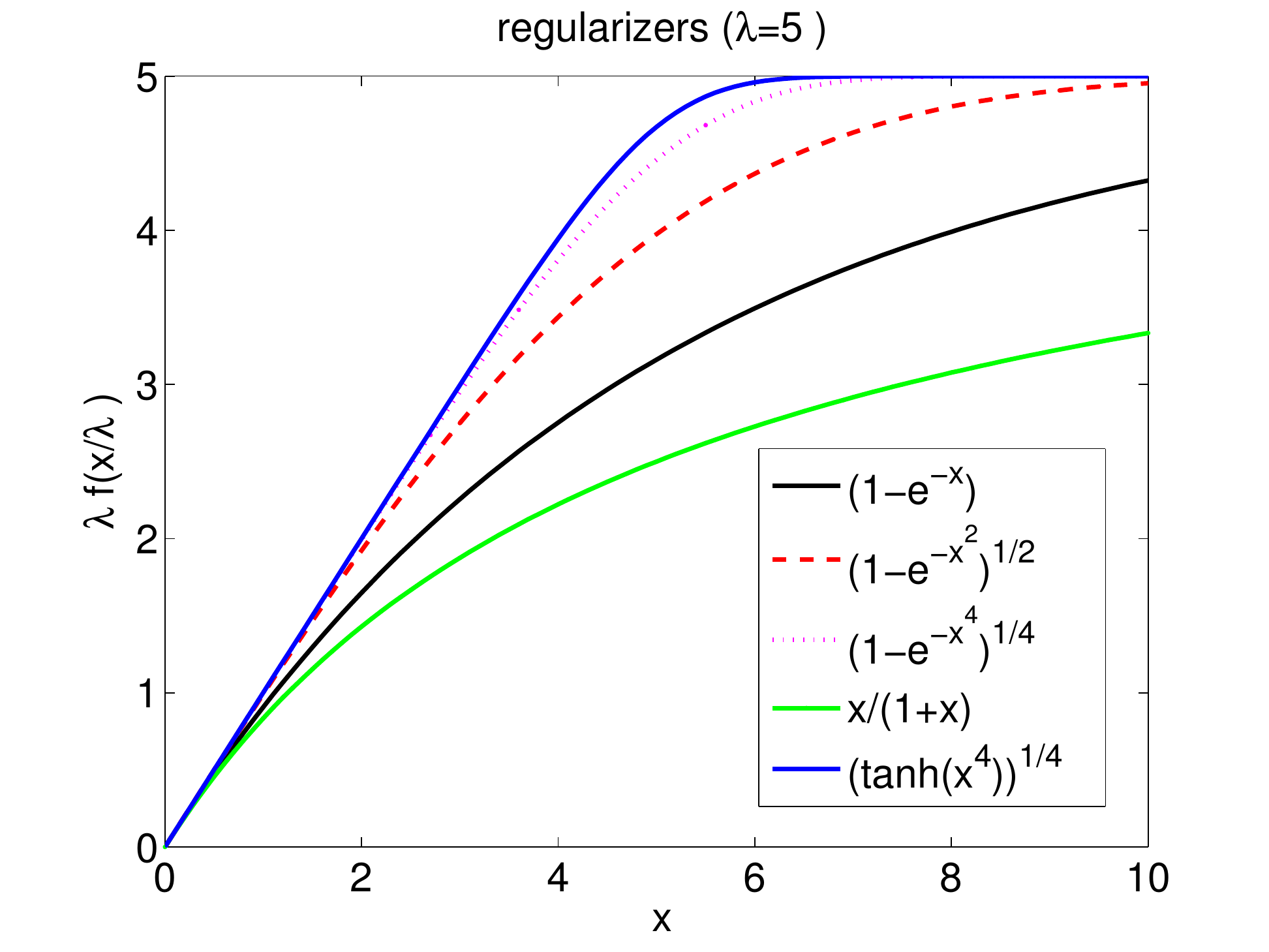}
\caption{Regularizers for eq. (\ref{eq:regularizer1}).}
\label{fig:regularizers}
\end{center}
\end{figure}

An obvious problem with deep holes, jumps or ragged noise in the
carrying capacity $\mathcal{K}$ is this: the right hand side of
equation (\ref{eq:godunovwithKC}) (as well as
(\ref{eq:godunovwithKTC}) in the next Section) contains the term
\begin{equation}
\left( 1 + \frac{k}{2} A_{y} + \frac{h}{2} (1-\frac{u^{\star}}{\mathcal{K}})
\right) u^{\star},
\label{eq:badeulerterm}
\end{equation}
which for large ${u^{\star}}/{\mathcal{K}}$ basically determines
the sign of the $u(t+h)$ on the left side of (\ref{eq:godunovwithKC})
(likewise (\ref{eq:godunovwithKTC})).  Since both $u^{\star}$ and
$\mathcal{K}$ are positive, if $\mathcal{K}$ is very small in some
pixel, then we may have \[ \frac{u^{\star}}{\mathcal{K}} \gg
\frac{2}{h},\] in which case $u(t+h)$ changes sign. This is
unphysical, so we would like to regularize the term
(\ref{eq:badeulerterm}).
To do so, we need to find a monotonically increasing function, call it
$g(x)$, such that
\[               
g(x) \sim \left\{ 
\begin{array}{ll}
x & \mbox{when} \, x \, \mbox{is small} \\
1 & \mbox{when} \, x \, \mbox{is large.} 
\end{array} 
\right. 
\] 
Multiple choices are available, as shown in
Figure~\ref{fig:regularizers}. We want a straight line with slope 1
when $x$ is small, then a smooth but rapid cut-off when 
$x = \frac{u^{\star}}{\mathcal{K}}$ gets too large. Some suitable choices
are $(1-\exp{(-x^{\beta})})^{1/\beta}$, or just as cheap to compute,
\begin{equation}
g(x) = \left( \tanh{(x^{\beta})} \right)^{1/\beta},
\label{eq:regularizer0}
\end{equation}
or any variant shown in Figure~\ref{fig:regularizers}. We choose
$g(x)$ with $\beta=4$. The regularization (\ref{eq:regularizer0}) to be
used in (\ref{eq:badeulerterm}) and thus (\ref{eq:godunovwithKC})
(likewise (\ref{eq:godunovwithKTC}) in the next section) is then
modified to
\begin{equation}
\frac{u^{\star}}{\mathcal{K}} \mapsto \frac{1}{h}
g(h \frac{u^{\star}}{\mathcal{K}})
\label{eq:regularizer1}
\end{equation}
where $h$ is the step-size. In the desert test presented above, using
the regularizer yields the same results as decreasing $h$ tenfold.

\begin{figure}
\begin{center}
\begin{minipage}{0.48\textwidth}
  \includegraphics[width=\textwidth]{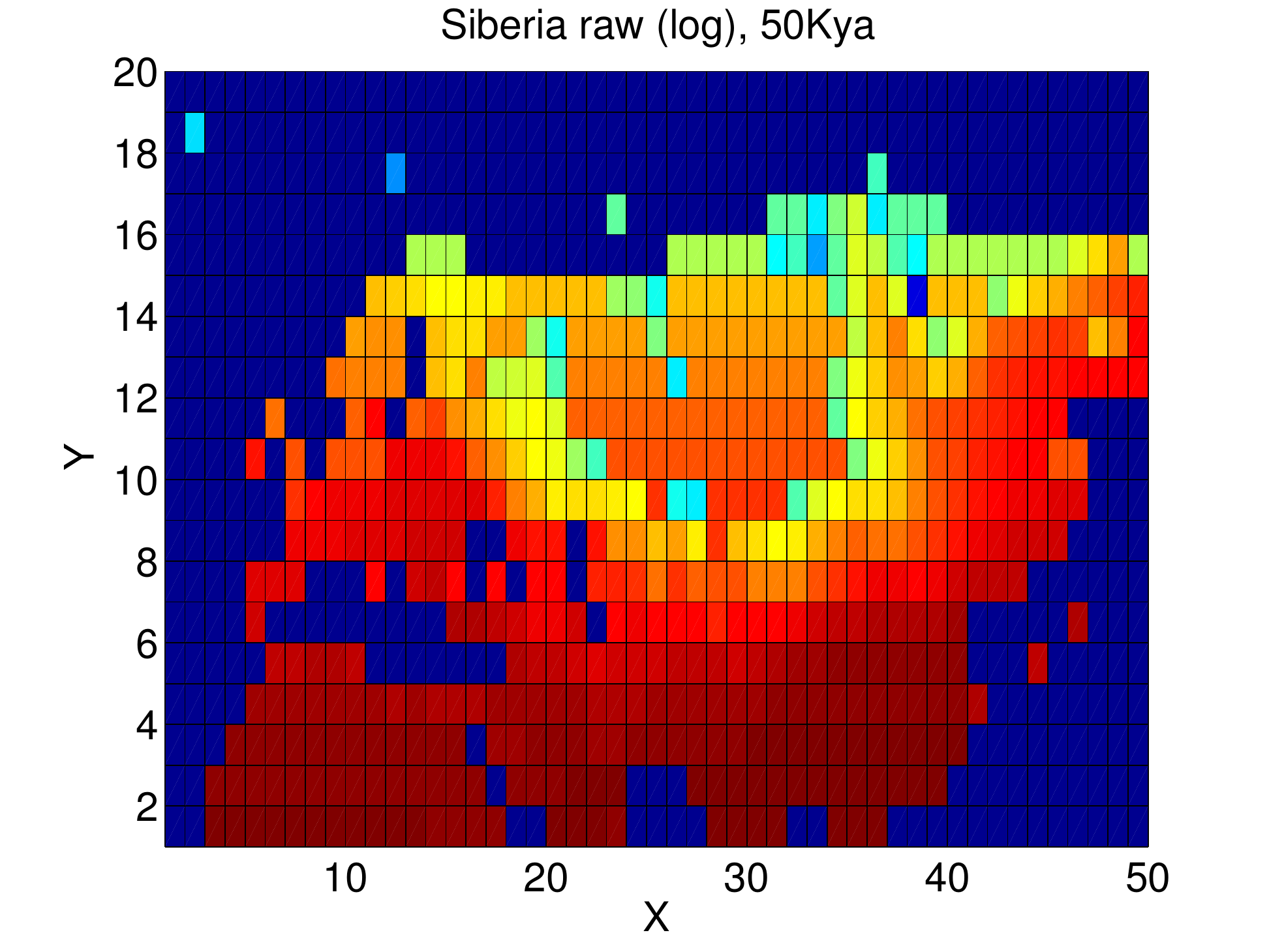}
\end{minipage}
\hfill
\begin{minipage}{0.48\textwidth}
  \includegraphics[width=\textwidth]{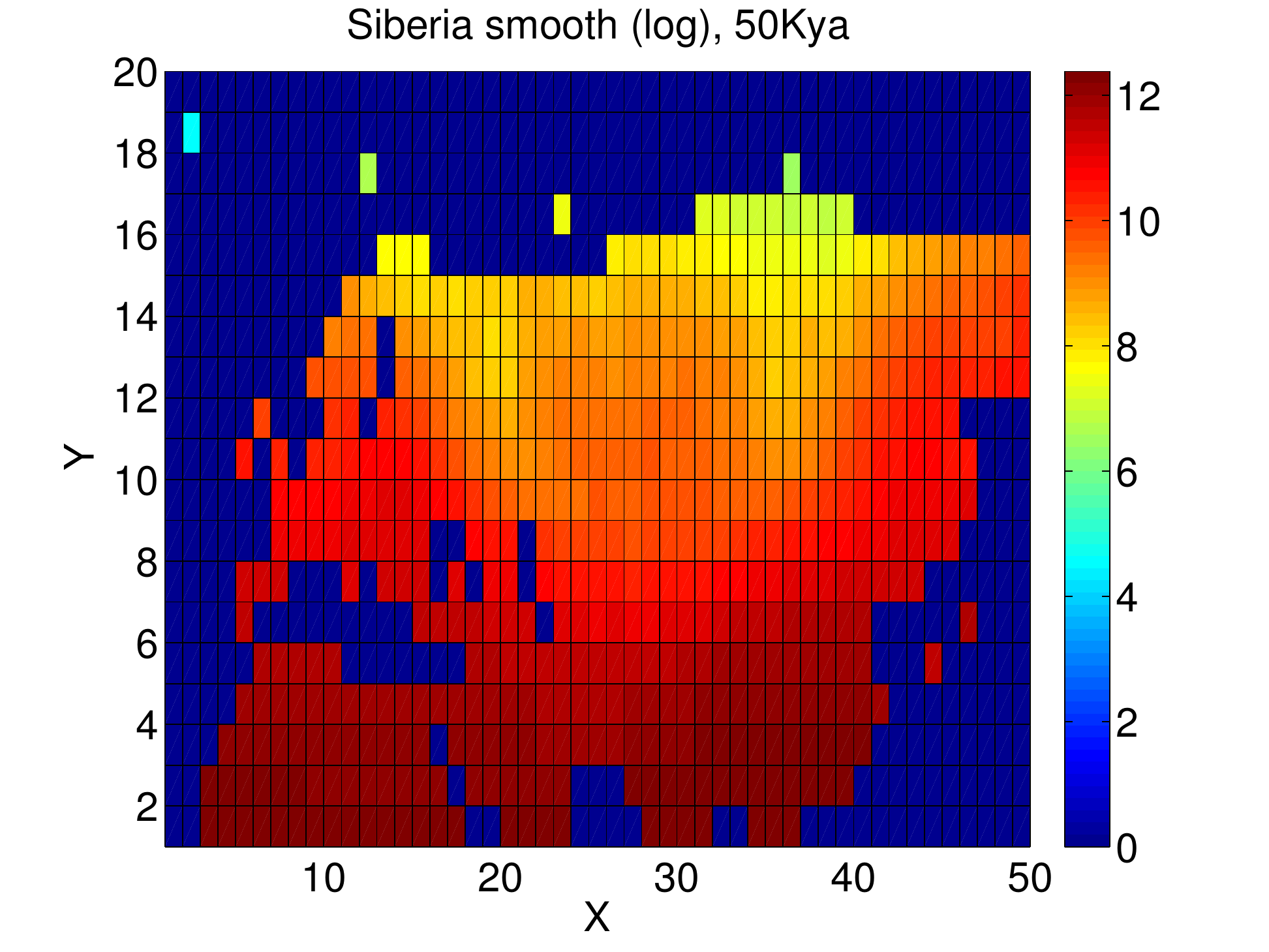}
\end{minipage}
   \caption{Left: Siberia rough 50kya NPP distribution. Note the
     $\mathcal K=0$ hole at pixel $\mbox{north}=14$, $\mbox{east}=38$.
     Right: same sub-map but smoothed by the low pass filter
     (\ref{eq:weightfunction}). Scales are
     $\log{(1+10\mathcal{K})}$ with the original $0.3 \leq \mathcal{K}
     \leq 3 \times 10^{4}$ biomass units (gm/$\mbox{m}^{2}$)
     \citep{grieser2006}. Both plots use the same colorbar.}
\label{fig:siberiaCCrawvssmooth}
\end{center}
\end{figure}

\begin{figure}
\begin{center}
\begin{minipage}{0.48\textwidth}
  \includegraphics[width=\textwidth]{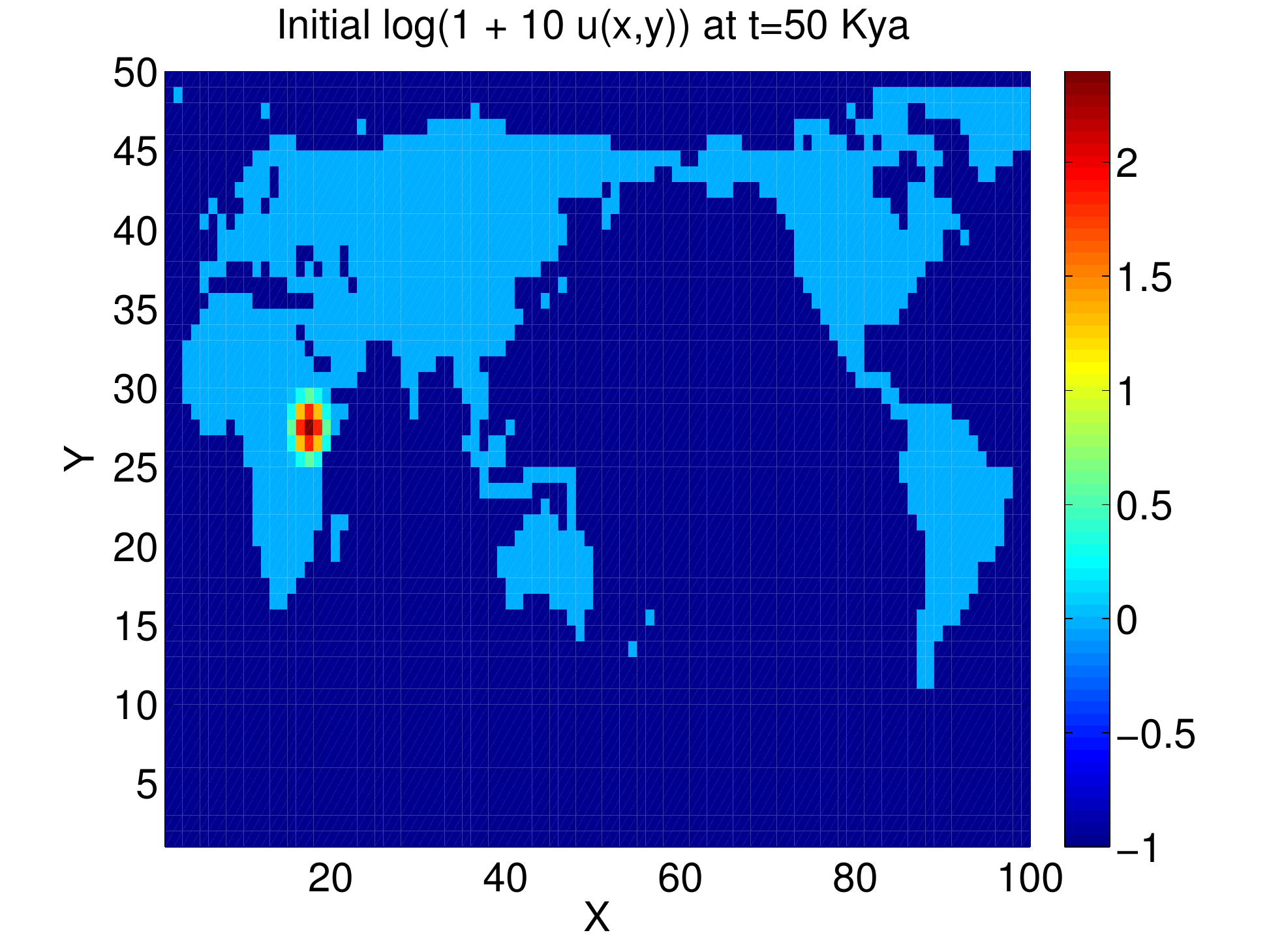}
\end{minipage}
\hfill
\begin{minipage}{0.48\textwidth}
  \includegraphics[width=\textwidth]{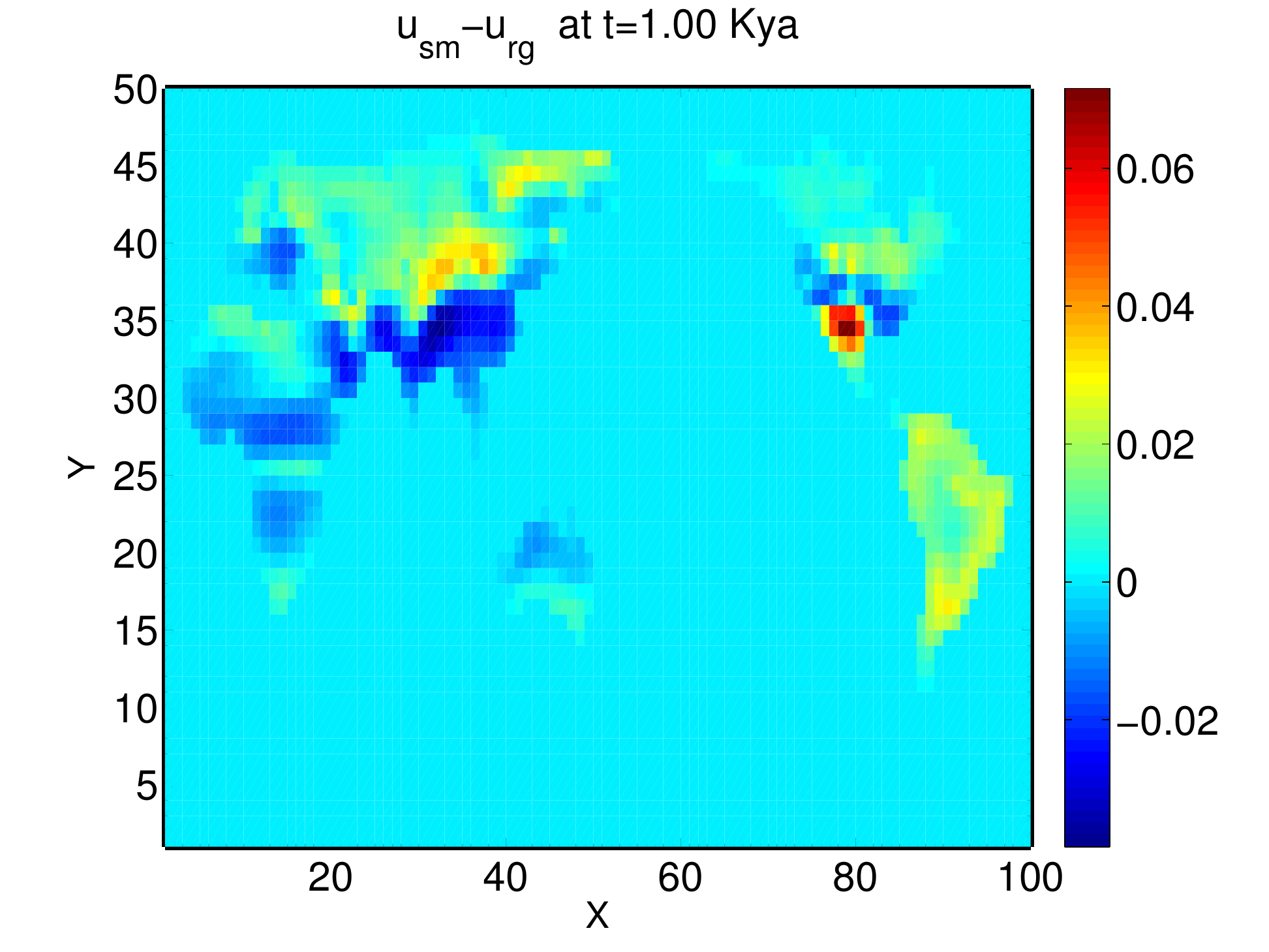}
\end{minipage}
\caption{Color-coded $\log{(1+10u)}$ plot for out-of-Africa dispersal.
  Dark blue shows {\em water}, with $\mathcal K=0$. Growth rate is
  $\lambda = 1.67\cdot10^{-3}\,\mathrm{ky}^{-1}$, and diffusion coefficient
  $c=208\,\mathrm{km}^2\mathrm{/ky}$ \citep{young1995}. Left: initial
  distribution at $t=50$~kya. This color scale is used in all
  subsequent dispersal plots below. Right:  differences between 
  the regularization (\ref{eq:regularizer0}) of $u$ 
 vs. smoothing of $\mathcal{K}(\mathbf{x},t)$ (\ref{eq:weightfunction}) at $t=1$~kya.}
  \label{fig:outofafrica0}
\end{center}
\end{figure}

One alternative to the above regularization (\ref{eq:regularizer1}) is to use 
smoothing by a low pass filter which weights a center pixel (map
coordinate $i,j$) and its nearest $(2L-1) \times (2L-1) - 1$
neighbors.  A neighboring pixel with $X,Y$ distances $i_{x}, j_{y}$
from $(i,j)$, has weight
\begin{equation}
w(i_{x},j_{y}) = \left( 1 - (i_{x}/L)^{2} \right)
\left( 1 - (j_{y}/L)^{2} \right)
\label{eq:weightfunction}
\end{equation}
for all $-L < i_{x} < L$ and $-L < j_{y} < L$, including center at
$i_{x} = j_{y} = 0$. Neighboring pixels with coordinates 
$(i+i_{x},j+ j_{x})$ having zero value (e.g. water), 
$i+i_{x} > {\rm N_X}$ or $i+i_{x} < 1$, 
are ignored.  In Mercator projection maps, $Y$ coordinates are
periodic in the longitude direction).  For each $(i,j)-$pixel to be
smoothed, a total of each accepted (non-zero) neighbor's weight was
kept and the resulting total was normalized appropriately. The choice
(\ref{eq:weightfunction}) is an approximate Gaussian weight 
$\exp{( -\Delta x^{2} - \Delta y^{2})}$, cut off at distances 
$|i_{x}| \geq L$ and $|j_{y}| \geq L$ or at uninhabitable
pixels.  Figure~\ref{fig:siberiaCCrawvssmooth} shows an example of the effect
of a low-pass filter (5 cells half-width, eq.~(\ref{eq:weightfunction})) 
smoothing on one of the maps used in this study. A comparison between
regularization and smoothing is shown on the right hand plot
of Figure~ \ref{fig:outofafrica0}.

\subsection{Time-dependent capacity maps}
\label{sec:fisherkppvarKxt}
Fluctuations in climate produce environmental changes in vegatation,
sea levels, opening/closing of land bridges, waxing/waning of ice sheets, and
perturbations to habitable areas in general. Thus, time-dependent
environments compel us to extend our procedure (\ref{eq:godunovwithK})
for both space- and time-dependent $\mathcal{K}(\mathbf{x},t)$ (see
Section~\ref{ss:timeinterpol}).

Since (\ref{eq:godunovwithK}) is basically the trapezoidal method (see
section 5.3 in \citep{leveque2007}), the modification for a
time-dependent $\mathcal K$ is as follows:
\begin{subequations} \label{eq:godunovwithKT}
\begin{align}
\left( 1 - \frac{k}{4} A_{x} \right) u^{\star} & = \left( 1 + \frac{k}{4} A_{x} 
\right) u(t) \label{eq:godunovwithKTA} \\
u_{E} &= u^{\star} + k A_{y} u^{\star} + h (1-\frac{u^{\star}}{\mathcal{K}(t)})
u^{\star} \label{eq:godunovwithKTB} \\
\left( 1 - \frac{k}{2} A_{y} - \frac{h}{2} (1-\frac{u_{E}}{\mathcal{K}(t+h)}) 
\right) u^{\star \star} & =
\left( 1 + \frac{k}{2} A_{y} + \frac{h}{2} (1-\frac{u^{\star}}{\mathcal{K}(t)}) 
\right) u^{\star} \label{eq:godunovwithKTC} \\
\left( 1 - \frac{k}{4} A_{x} \right) u(t+h) & = 
\left( 1 + \frac{k}{4} A_{x} \right) u^{\star \star},
\label{eq:godunovwithKTD}
\end{align}
\end{subequations}
where we have again suppressed the $\mathbf{x}$ dependence of
$\mathcal{K}(\mathbf{x},t)$ for simplicity of notation.

\subsection{Fisher/KPP on geographical maps}
\label{subsec:geographicalmaps}
In order to solve Fisher/KPP on an irregular domain such as a
geographical map, it is sufficient to break down the map into $x$- and
$y$-direction segments, imposing a $u=0$ boundary condition at their
endpoints; the solver can act on each segment independently,
alternating the direction \citep{leveque2007,ritchmyer1967}
of integration as discussed above. This
approach lends itself also to efficient
parallelization. Appendix~\ref{sec:code} and
Figure~\ref{fig:mapsegmentation} illustrate in detail our scheme with
a sample Matlab code.

Note that holes in $\mathcal K$ maps can represent the partial closing
of land bridges without the necessity of re-segmenting land portions,
as in Section~\ref{sec:fisherkpp} and shown in
Figure~\ref{fig:mapsegmentation}. However, the regularization scheme
(Section~\ref{sec:noisyNPPmaps}) should be used with caution. A $u=0$
boundary condition region is not the same as one with low
$\mathcal{K}$, which can pass a tiny population into a subequent region
where it may well flourish. For example, historically known falling
sea levels opened passages across Bab al Mandeb to open South Asia for
human dispersal \citep{siddall2003}, and the retreat of the North
American ice sheet opened a passage on the Bering strait around 12~kya
(kilo-years ago).

\section{World-wide hominin dispersal}
\label{sec:dispersal}
We now turn to a sample application of the methods presented above:
the world-wide dispersal of \emph{Homo sapiens} out of
Africa.

\subsection{Capacity maps}

Our construction of a time-dependent $\mathcal{K}$ uses Net Primary
Productivity (NPP) as a proxy \citep{eriksson2012}. The Miami model
\citep{grieser2006} was originally formulated in 1972 to estimate NPP
(in grams (of carbon) in dry organic matter/m$^2$/day) from annual 
temperature and rainfall \citep{lieth1975}. In order to compute our 
NPP maps, we obtained the temperature and precipitation data 
from simulations by the {\sc bridge} program \citep{bridgeproj}
organized at the University of Bristol
% \footnote{\url{http://www.bridge.bris.ac.uk}}
\citep{bigelow2003, harrison2001, pickett2004, prentice2000}. The
simulation data that we used were computed on a $96 \times 73$ grid,
which we interpolated to size $100 \times 50$ and converted to NPP
maps by applying the formulas given in \citep{grieser2006}.  In
Figure~\ref{fig:siberiaCCrawvssmooth},
the original NPP units in (grams of C)/m$^{2}/day$ were adopted.

World-wide NPP data are difficult to obtain, so our Miami model-like
maps are rough. As we showed in Section~\ref{sec:gaptest}, our Godunov
solver is fairly robust with respect to abrupt steps in the carrying
capacity $\mathcal{K}(\mathbf{x})$. However, an alternative to the
regularization scheme proposed in Section~\ref{sec:noisyNPPmaps} would
be desirable, when dealing with maps in which noise and holes may not
necessarily correspond to real physical conditions of the model
system. 

\subsection{Time interpolation of maps}\label{ss:timeinterpol}
We assembled 61 NPP maps, from 120 kya to 1 kya. These are in 4 ky
steps for the first 10, 2 ky steps for the next 21, then 1 ky for the
remainder. Since the time stepper in our Godunov scheme can have no
information about the continuity in time of NPP maps, an interpolation
scheme needs to be used.

If some estimates or proxies for the carrying capacities 
$\mathcal{K}_{L}$, $\mathcal{K}_{H}$ at times $t_L$, $t_H$ are available, one possible
estimate for $\mathcal{K}(t)$ is a homotopy
\[
 \mathcal{K}(\mathbf x, t) = \mathcal{K}_{L}(\mathbf x)(1-S(t)) 
+ \mathcal{K}_{H}(\mathbf x) S(t)
\] 
where a sigmoid function $0 = S(t_L) <  S(t) < S(t_H) = 1$ will 
smoothly interpolate between the two time
frames. There are many choices available, such as that used in
\citep{eriksson2012}. For this study, we used the following variant.

\begin{figure}
\begin{center}
  \includegraphics[width=0.7\textwidth]{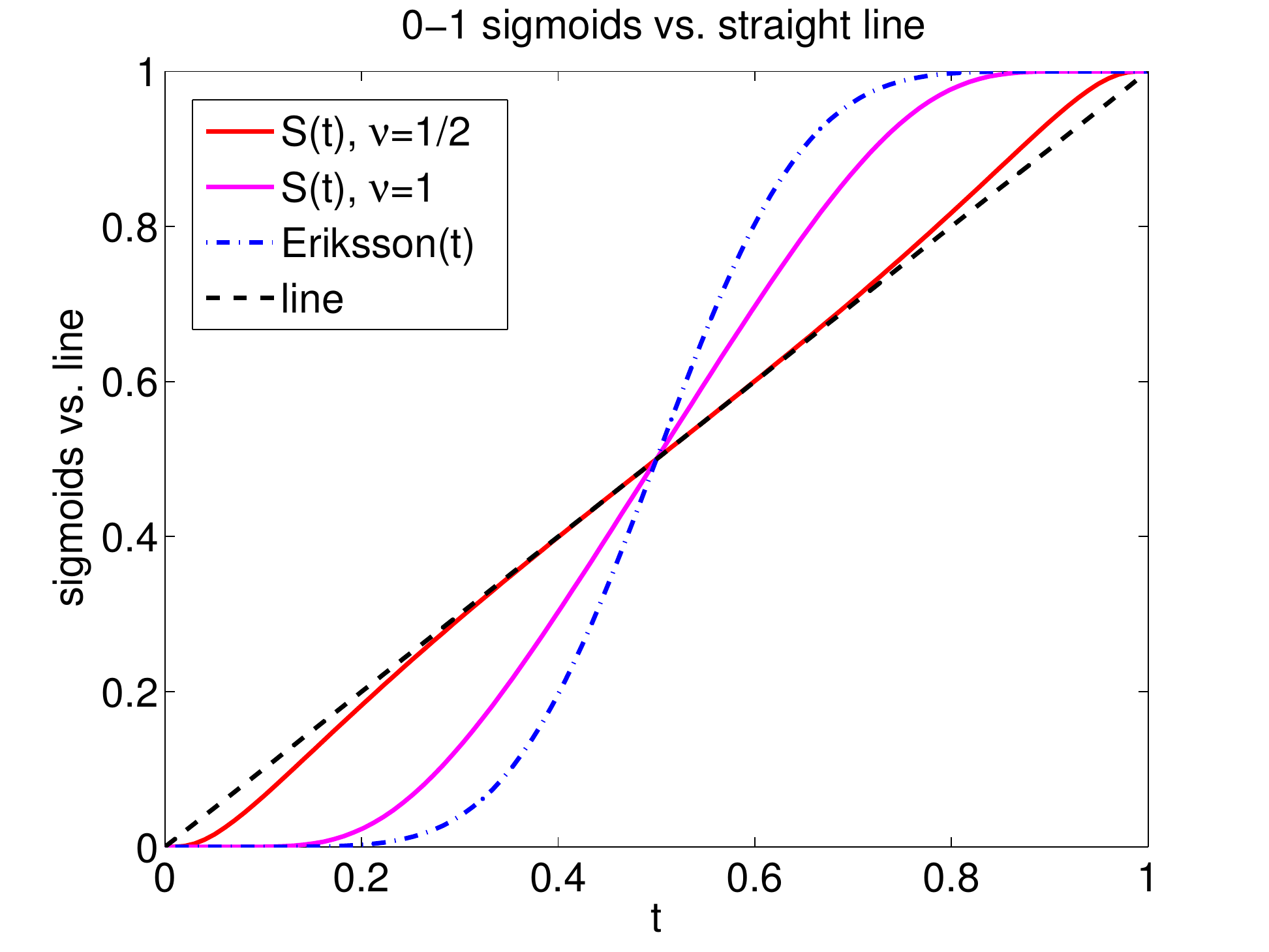}
\caption{Interpolation between $\mathcal K$ maps at different
  times: sigmoids (\ref{eq:ourinterpolator}) for $\nu=1/2$ and
  $\nu=1$, straight line, and Eriksson's $f(f(f(t)))$
  model \citep{eriksson2012}.}
\label{fig:sigmoids}
\end{center}
\end{figure}

Start with the classical sigmoid
\begin{equation}
S(z) = \left( 1 + e^{-z} \right)^{-1}
\label{eq:classicalsigmoid}
\end{equation}
which is zero at $z = -\infty$ and unity at $z=+\infty$. The $-\infty
< z < +\infty$ interval is not what we want, but with a small change
the following $t \mapsto z$ transformation permits many variants:
\begin{align}
z & = \frac{2 \Delta T (t-t_{L}) -\Delta T^2}{((t-t_{L})(\Delta T - 
(t-t_{L})))^\nu}.
\label{eq:ourinterpolator}
\end{align}
where $\Delta T = t_{H} - t_{L}$. Notice that $S(z(t_{L})) = 0$ and
$S(z(t_{H})) = 1$. The exponent $\nu$ in (\ref{eq:ourinterpolator})
gives some freedom in choosing a particular form for $z$ for almost
any $\nu > 0$. If $\nu < 1/2$, ${\rm d}^{2} S/{\rm d}t^{2}$ has more
than two sign changes, so $\nu \geq 1/2$ is preferable.  By the choice
$\nu=1/2$, the interpolant is nearly a straight line: see
Figure~\ref{fig:sigmoids}. However, its turn-up at $t=t_{L}$ and
turn-down at $t=t_{H}$ numerically resemble very quick derivative
changes. Thus, for this test case we choose $\nu=1$.  At both ends,
all derivatives of $S(z(t))$ in $t$ vanish.  Also notice the
forward/backward symmetry $S(z(t_{H}-t)) = 1 - S(z(t))$ for $t_{L}
\leq t \leq t_{H}$ \citep{eriksson2012}.

\subsection{Population dispersal}
\begin{figure}
\begin{center}
\begin{minipage}{0.48\textwidth}
  \includegraphics[width=\textwidth]{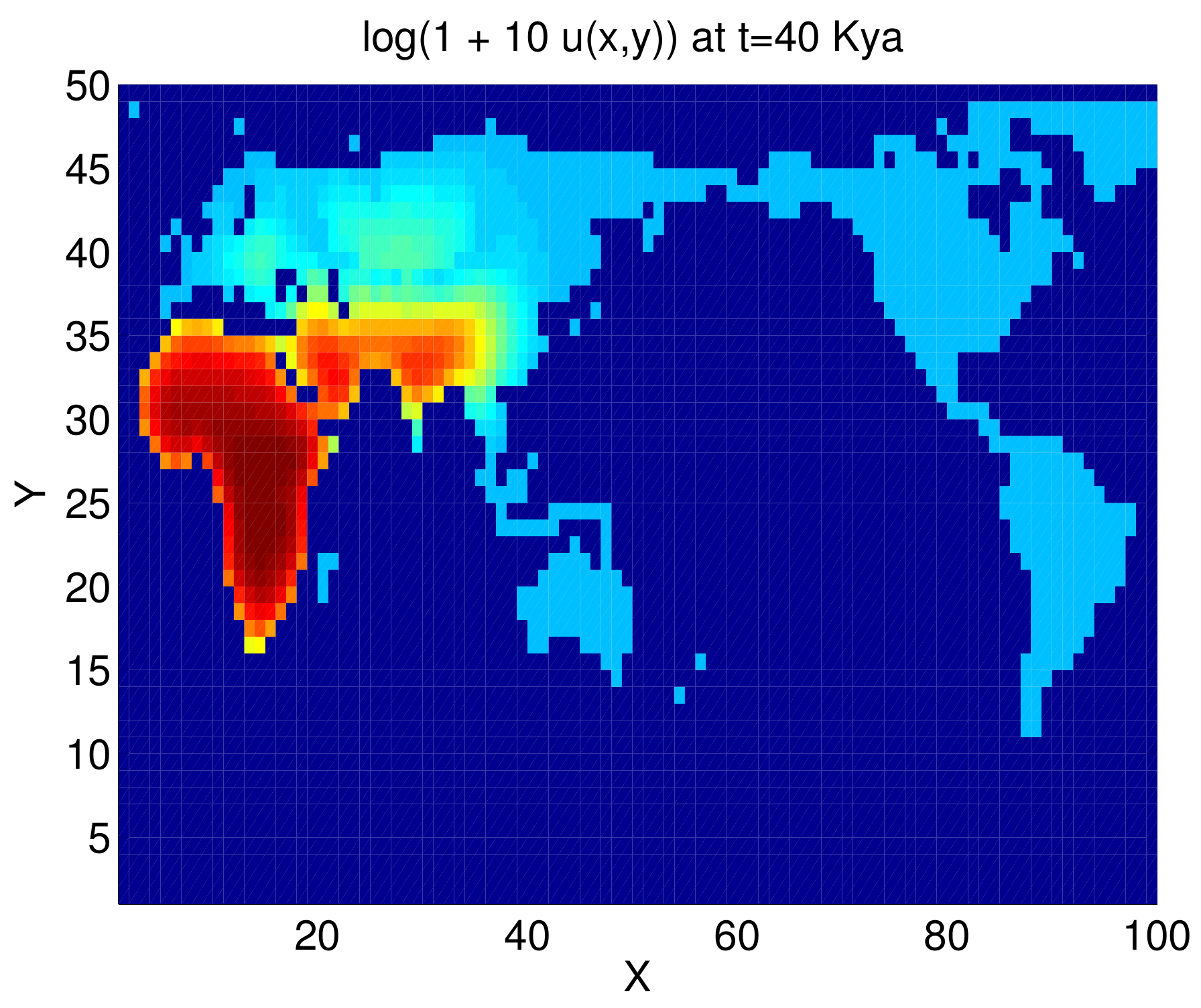}
\end{minipage}
\hfill
\begin{minipage}{0.48\textwidth}
  \includegraphics[width=\textwidth]{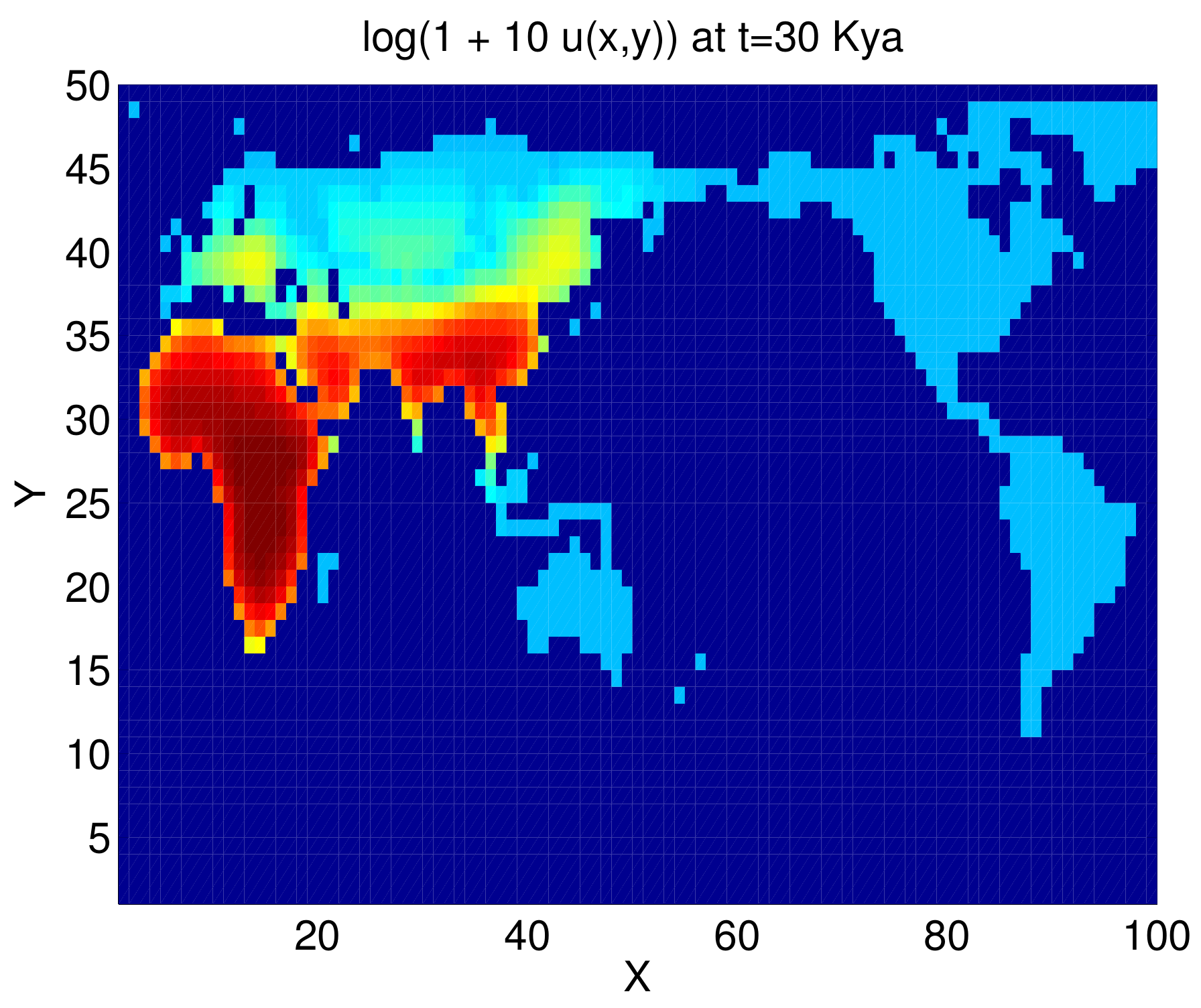}
\end{minipage}
\caption{Out of Africa simulations by regularization 
  (\ref{eq:regularizer0}). Left: 
  population at $t=40$~kya. Right: population at $t=30$~kya.
  Color scale is from Fig.~\ref{fig:outofafrica0} (left).} 
  \label{fig:outofafrica1}
\end{center}
\end{figure}

\begin{figure}
\begin{center}
\begin{minipage}{0.48\textwidth}
  \includegraphics[width=\textwidth]{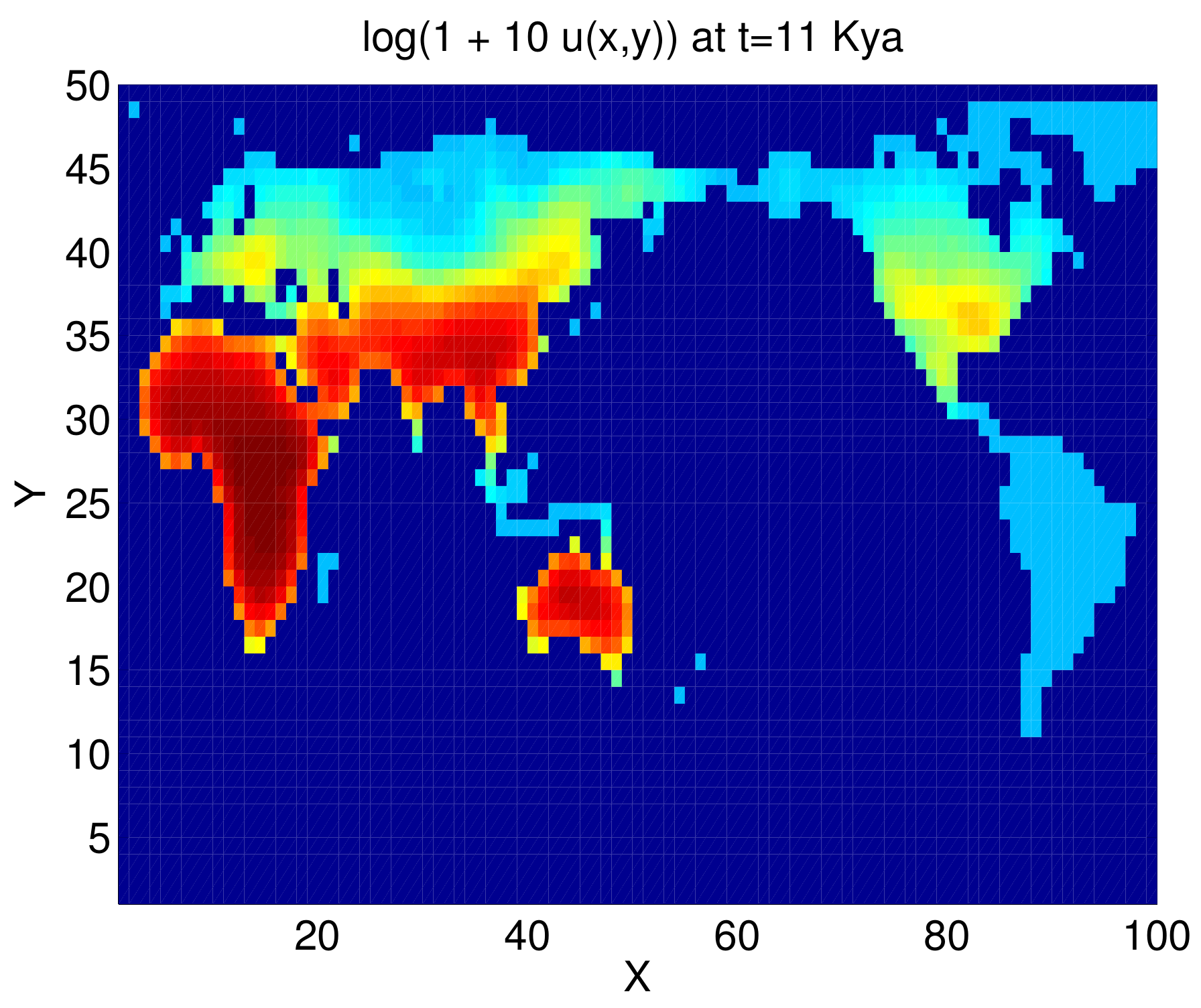}
\end{minipage}
\hfill
\begin{minipage}{0.48\textwidth}
  \includegraphics[width=\textwidth]{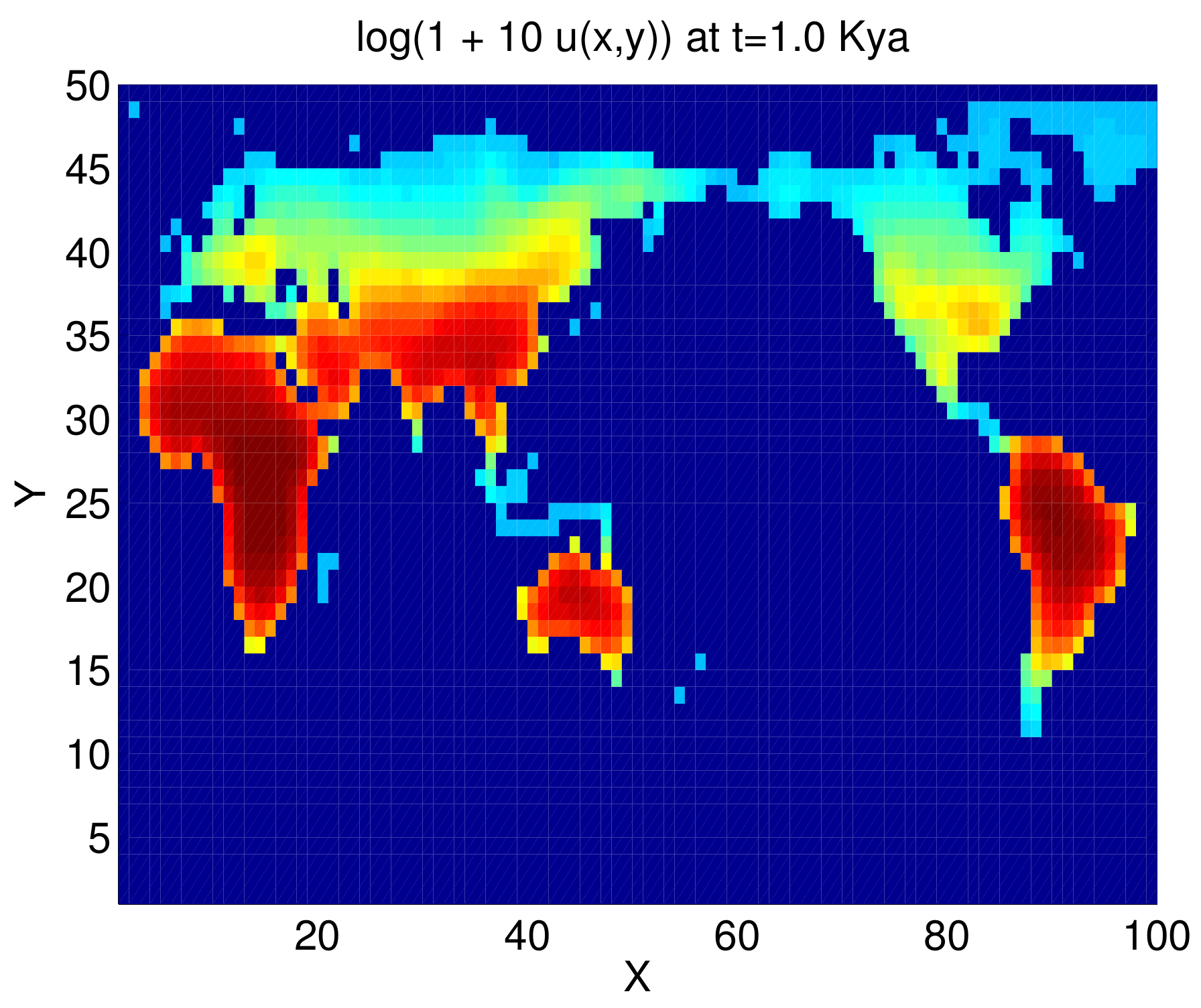}
\end{minipage}
\caption{Out of Africa simulations with regularization
  (\ref{eq:regularizer0}). Left: 
  population at $t=11$~kya. Right: population at $t=1$~kya.
  Color scale is from Fig.~\ref{fig:outofafrica0} (left).}
  \label{fig:outofafrica2}
\end{center}
\end{figure}

We set up the initial population density at $t=50$~kya as shown in the
top left panel of Figure~\ref{fig:outofafrica0}. The integration units
are scaled following Tab.~\ref{tab:scalings} such that $\lambda =
1.67\cdot10^{-3}$~ky$^{-1}$ and $c=208$~km$^2$/ky (consistent with those
used in \citep{young1995}). The solver is then run using time frame
$\mathcal K$ maps described above, down to 1~kya.  The remaining
plots (Figures~\ref{fig:outofafrica1},\ref{fig:outofafrica2}) 
display the resulting population dispersal simulation on 
unsmoothed $\mathcal{K}$ maps, regularized by
(\ref{eq:regularizer0}). Using population parameters consistent with
the literature, the gross features of the late (50-60 kya)
out-of-Africa dispersal of
{\it Homo sapiens} are reproduced \citep{forster2004}, e.g.~the
colonization of Western Europe by $\sim 40$~kya and that of South
America before 14~kya. Using the solver with the same initial
conditions but on a smooth NPP map, like the one shown in right panel
of Figure~\ref{fig:siberiaCCrawvssmooth}, {\it without} the regularizer
(\ref{eq:regularizer0}), yields the same wavefront propagation speed.

\section{Conclusions}
\begin{figure}
\begin{center}
\begin{minipage}{0.48\textwidth}
  \includegraphics[width=\textwidth]{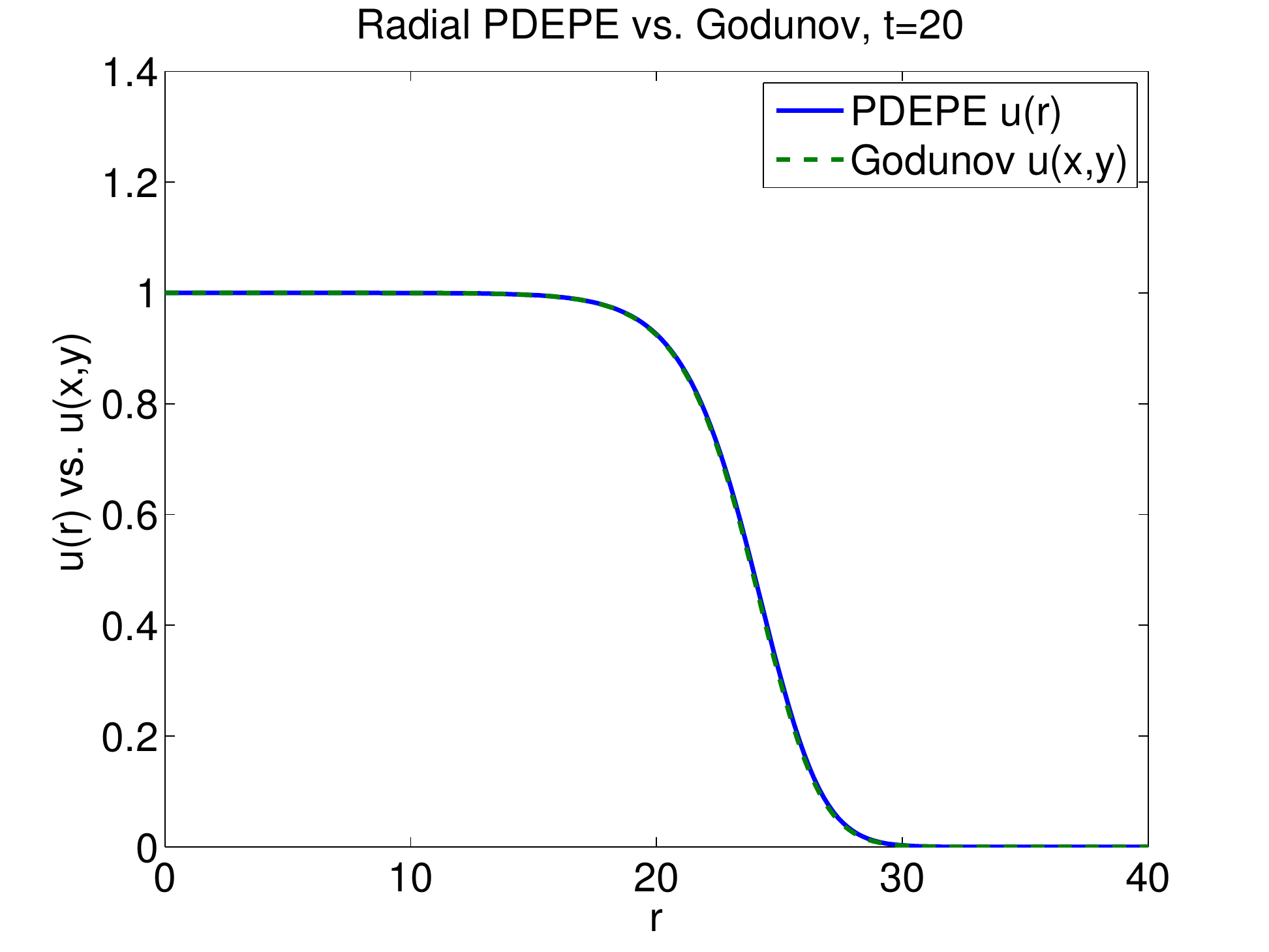}
\end{minipage}
\hfill
\begin{minipage}{0.48\textwidth}
  \includegraphics[width=\textwidth]{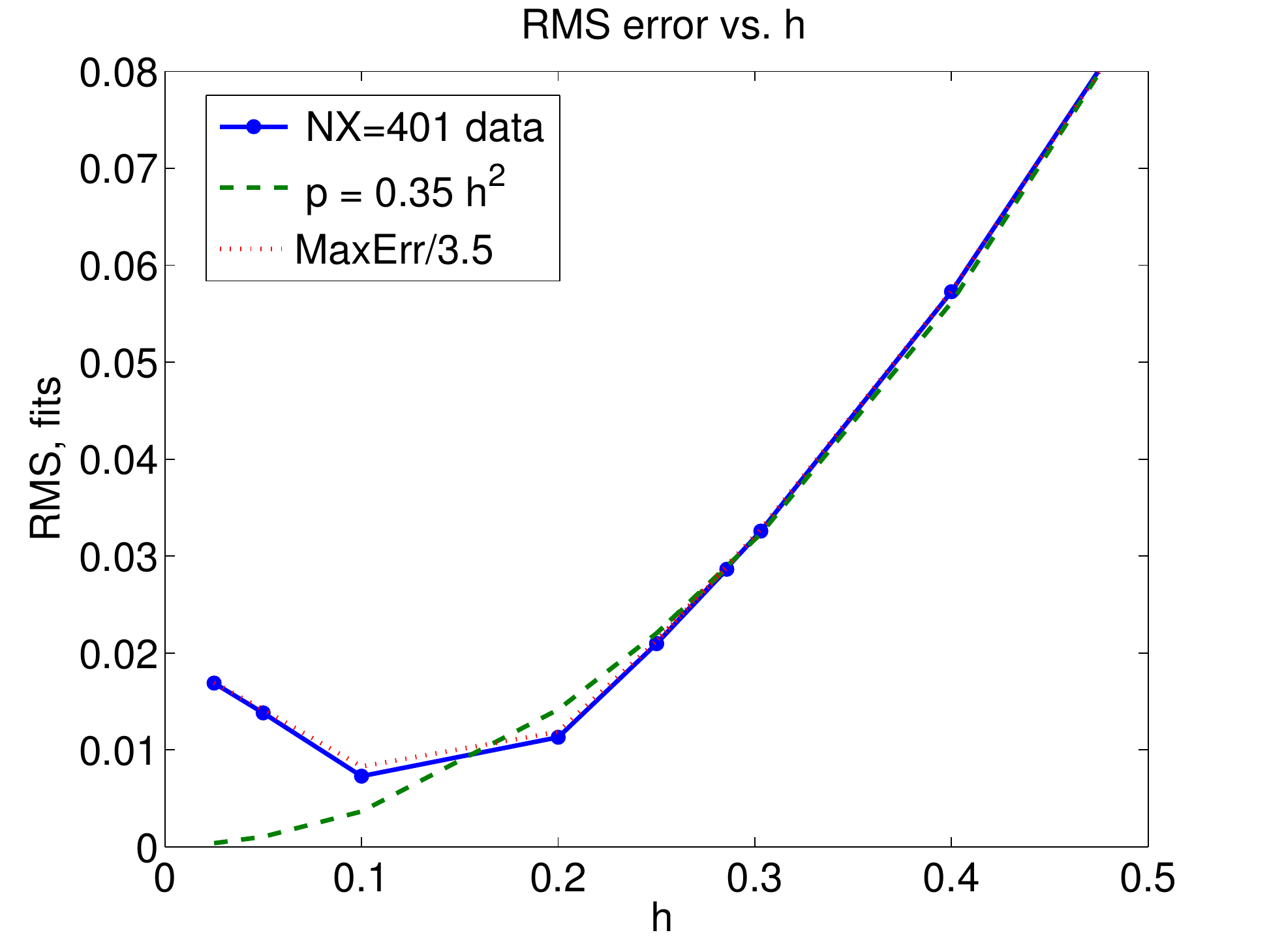}
\end{minipage}
\caption{Left: 2-D Wave front at $t=20$ compared to \texttt{pdepe}.
  Right: Rough RMS error as a function of step size, $h$, at $t=20$.
  Errors $\epsilon_{RMS}$ and $\epsilon_{MAX} \approx 3.5
  \epsilon_{RMS}$ have approximately $O(h^2)$ behavior for $h \geq
  1/10$, but increase if $h$ is too small: overall estimate is
  $O(h^{2}) + O(\Delta x^{2}) + O(\frac{1}{h} \Delta x^{3})$.}
\label{fig:2Dwavefrontanderrors}
\end{center}
\end{figure}

In this paper, we presented a novel semi-implicit Godunov scheme for
the Fisher/KPP equation with a constant carrying capacity
$\mathcal{K}$, described in Section~\ref{sec:methodsandsplittings}. In
one dimension, the expected traveling wave
\citep{fisher1937,kolmogorov1937} develops as shown in
Figure~\ref{fig:1dplots}.  In other tests, not shown here, we saw that almost
any concentrated inital condition will develop similar waves: for
example, two nearby peaks. Our
scheme is on a rectangular grid, so in 2-D we need to ensure that in
cylindrically symmetric situations we can control the errors due to
the $x - y$ asymmetry.  The error plots,
Figure~\ref{fig:symmetryornot}, show that we can reduce the
asymmetries somewhat by an alternating direction scheme. In any case,
these errors are very small even for a relatively high CFL number. 
As in the one-dimensional situation, a
traveling wave also develops as expected
\citep{kolmogorov1937}: Figure~\ref{fig:2dpcolorplots} and
Figure~\ref{fig:2Dwavefrontanderrors} show this development and
compare the results to the 2-D symmetric version of the Matlab
function \texttt{pdepe}. Because the $\mathbf{x}$ discretization has
truncation errors proportional to $(\Delta x)^{3}$, after
$O(\frac{1}{h})$ time steps we should not be surprised to see the
errors shown in Figure~\ref{fig:2Dwavefrontanderrors} turn up for very
small $h$ and behave roughly as $O(\frac{1}{h} (\Delta x)^{3})$.

In Secs.~\ref{sec:fisherkppvarK} we extended our procedure to handle
an $\mathbf{x}$-dependent $\mathcal{K}$, specifically
eq.~(\ref{eq:godunovwithK}).  In order to regularize our solver
against bad behavior when dealing with $\mathcal{K}$ maps inferred
from real-world data, in Section~\ref{sec:gaptest} we studied both
regularization and the expected constant velocity of the traveling
wave. Except for small transition regions entering or leaving a region
of low carrying capacity, the velocity ($\sqrt{2}$ in our scaling) is
indeed constant. In Section~\ref{sec:fisherkppvarKxt}, we
went further to develop a scheme for the situation with both $\mathbf{x}$
and $t$ dependent carrying capacity $\mathcal{K}$ via
eq.(\ref{eq:godunovwithKT}).

Finally, this scheme has been applied to a prototypical case in
population dynamics: the out-of-Africa dispersal of {\it Homo
  sapiens}. On the Mercator projected world map, by using vegetation
net primal productivity as a proxy for carrying capacity,
Figure~\ref{fig:outofafrica0} shows that by the regularization of
$\mathcal{K}$ in space via eq.~(\ref{eq:regularizer0}) but
interpolating in time yields stable and reasonable results. In fact,
the results showing ancestor arrival in NE India at roughly 40~kya, then crossing the 
Bering Strait before 10~kya, and multiple routes into South Asia \citep{reyescentento2014}
are very encouraging. Honesty requires that we admit our size (408km)$^{2}$
pixels do not resolve the two crossing points at Bab al Mendab and
Sinai adequately. Additionally, a {\em switch} 
(see: Section~\ref{subsec:geographicalmaps})
which would allow passage at the Bering Strait seems unnecessary due to the
interesting coincidence that the hominin wave front reaches this passage
at the beginning of the last ice age. If it were blocked previously, this
would have had no effect.

The core computations performed by our solver are \emph{independent}
tridiagonal solutions, which can be easily parallelized to deal with
larger grids. In order to further improve numerical performance, in
the Appendix~\ref{sec:code}, we discussed a compressed storage scheme to
integrate the Fisher/KPP equation on a projected world map (or any
other irregularly-shaped domain). The alternating direction scheme
discussed in Section~\ref{sec:Godunov2Dsplitting} also works with this
compressed storage. In the case of world dispersal discussed in
Section~\ref{sec:dispersal}, since about 71\% of the earth's surface is
water, this compressed storage reduces computational work by the same
amount.

\section*{Provenance}
For this paper, the simulations were run on either a Mac Mini, 2.4 GHz
Intel Core Duo (Mac OS 10.6.8), or a \texttt{MacBook} laptop with the 
same processor specifications but running OS 10.8.4.
On the Mini, \texttt{MatLab} (7.10.0.499) R2010a was used, respectively 
(8.1.0.604) R2013a on the MacBook laptop.

\section*{Acknowledgments}
The authors are thankful for the encouragement, financial
arrangements, and technical help from Prof.~G.~Lake, Prof.~S.~Mishra,
and Dr.~I.~Schloenvogt. Dr.~Schloenvogt's assistance at
critical junctures was particularly valuable.

\appendix
\section{Map segmentation appendix}\label{sec:code}
\begin{figure}
\begin{center}
\includegraphics[width=0.7\textwidth]{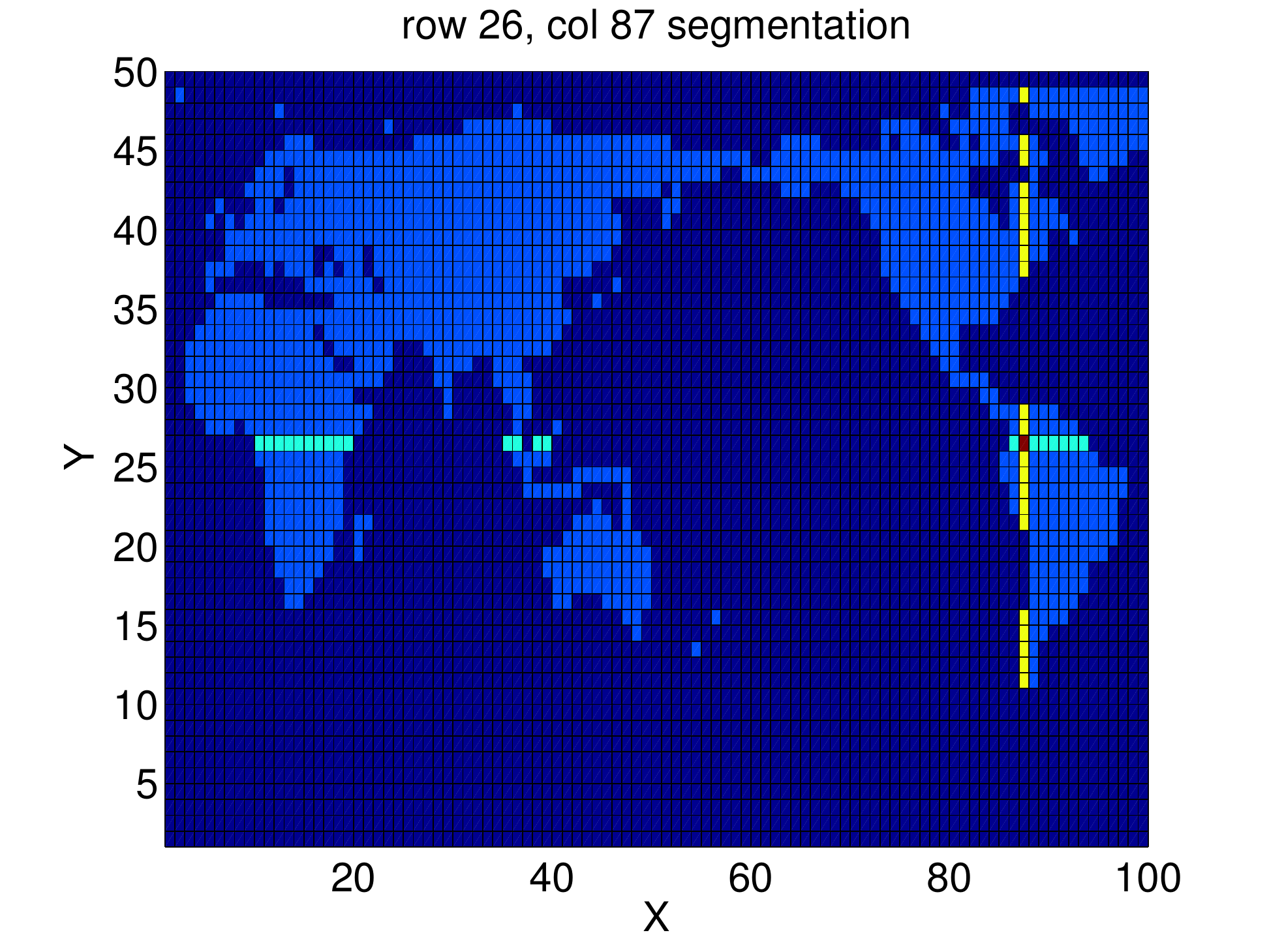}
\caption{Map with habitable regions set to 1 (light blue) and water
  regions set to 0 (dark blue) on a $N_x=100$, $N_y=50$ grid. The figure
  shows the segmentation for row 26 and column 87, with 4 and 5
  segments, shown in cyan and yellow respectively.}
\label{fig:mapsegmentation}
\end{center}
\end{figure}

In order to use our solver on a geographical map, it suffices to use a
map outline, i.e., a rectangular grid with 1's in habitable regions,
and 0's in the water, as in Figure~\ref{fig:mapsegmentation}. Since
each Godunov direction step only involves a row, or a column,
independently, we can set up the following indexing scheme.  Each row
$i = 1 \dots 50$ in Figure~\ref{fig:mapsegmentation} will have
\texttt{nysegs(i)} of habitable segments, whose starting and ending
positions are \texttt{ystart\_seg(k)} and \texttt{yend\_seg(k)},
respectively, where $k=1 \dots$\texttt{nysegs(i)}.  Likewise, for $j=1
\dots 100$ columns each with \texttt{nxsegs(j)} also with start and
end positions. Roughly 24\% of the world map is land, i.e. habitable.
As an example, Figure~\ref{fig:mapsegmentation} shows row 26 has 4
segments of varying size. Likewise, column 87 has 5 segments. A sample
\texttt{MatLab} code in Appendix~\ref{sec:code} illustrates the scheme
for 1/2 step of $x$-direction updates, followed by a full step of
$y$-direction updates, then again 1/2 step of $x$-direction segment
updates. The alternating direction method, described in
Section~\ref{sec:Godunov2Dsplitting}, just interchanges $x
\leftrightarrow y$ on alternate time steps.

For the reader's convenience, we include here a sample Matlab code
of our Godunov-Strang-Yoshida scheme.
\begin{verbatim}

% NY X direction updates for half-step1
      locx = 0;
      for j=1:NY
          nsegs = nxsegs(j);
          for k=1:nsegs
              istart=xstart_seg(locx+k); iend=xend_seg(locx+k);
              ninseg=iend-istart+1;
              u0(1) = 0; u0(ninseg+2) = 0;  % boundary values 
              u0(2:ninseg+1) = u(istart:iend,j);
% eq. (15a) solution:
              ut = godunovstep1(ninseg+2,h,kcfl,u0,sc1,sc2);
              u(istart:iend,j) = ut(2:ninseg+1);
          end
          locx = locx + nsegs;
      end
% NX Y direction updates step2
      locy = 0;
      for i=1:NX
          nsegs = nysegs(i);
          for k=1:nsegs
              jstart=ystart_seg(locy+k); jend=yend_seg(locy+k);
              ninseg=jend-jstart+1;
              u1(1) = 0; u1(ninseg+2) = 0; % boundary values
              u1(2:ninseg+1) = u(i,jstart:jend)';
              cap(1) = 1; cap(ninseg+2) = 1;
              cap(2:ninseg+1) = kap(i,jstart:jend);
% eq. (15c) solution:
              ut = godunovstep2(ninseg+2,h,kcfl,u1,sc1,sc2,cap);
              u(i,jstart:jend) = ut(2:ninseg+1);
          end
          locy = locy + nsegs;
      end
% Repeat godunovstep1, as above for eq. (15d)
%     locx = 0;
%     for j=1:NY
%         ETC 
%     end
\end{verbatim}

\pagebreak
\bibliographystyle{plainnat}
\bibliography{splitting}

\pagebreak

\end{document}